\documentclass[11pt, reqno]{amsart}

\evensidemargin
\oddsidemargin
\usepackage{amsmath, amsthm, epsfig, color, latexsym, amssymb, bbm}
\usepackage{xspace, tensor}
\usepackage{version} 
\usepackage{graphicx}

\newcommand{\assign}{:=}

\newcommand{\tmtextsf}[1]{\mathrm{#1}}
\newcommand{\wep}{\emph{wep}\xspace}
\newcommand{\whp}{\emph{whp}\xspace}

\newcommand{\eps}{\varepsilon}
\newcommand{\ssup}[1] {{\scriptscriptstyle{({#1}})}}
\newcommand{\one}{{\mathbbm 1}}

\newcommand{\R}{\mathbb R}
\newcommand{\B}{\mathcal B}

\newcommand{\skrin}{\mathcal N}

\newcommand{\E}{\mathbb E}
\renewcommand{\P}{\mathbb P}
\newcommand{\tmP}{\mathbb{P}}

\renewcommand{\phi}{\varphi}

\newcommand{\F}{\mathfrak F} 

\newcommand{\N}{\mathbb N}
\renewcommand{\P}{\mathbb P}
\newcommand {\T}{\mathbb T} 

\newcommand{\X}{\mathcal X}
\newcommand{\Y}{\mathcal Y}    

\newcommand{\V}{\mathcal V}
\newcommand{\W}{\mathcal W}
\newcommand{\PPP}{Poisson point process\xspace}
\newcommand{\indeg}{Z}

\newcommand{\x}{\mathbf{x}}
\newcommand{\y}{\mathbf{y}}
\newcommand{\z}{\mathbf{z}}

\newcommand{\Pp}{\P_*}   

\newtheorem{theorem}{Theorem}
\newtheorem{rem}[theorem]{Remark}

\newtheorem{corollary}[theorem]{Corollary}
\newtheorem{lemma}[theorem]{Lemma}
\newtheorem{defin}[theorem]{Definition}

\newtheorem{prop}[theorem]{Proposition}
\newtheorem{proposition}[theorem]{Proposition}

\theoremstyle{remark}

\newcommand{\heap}[2]  {\genfrac{}{}{0pt}{}{#1}{#2}}
\newcommand{\sfrac}[2] {\mbox{$\frac{#1}{#2}$}}

\begin{document} 

\ \\[-0.5cm]

\title[Robustness of scale-free spatial networks]
{Robustness of scale-free spatial networks}
\ \\[-1.5cm]

\author[Emmanuel Jacob and Peter M\"orters]{Emmanuel Jacob and Peter M\"orters}

\maketitle

\vspace{-0.6cm}



\begin{quote}
{\small {\bf Abstract:} }
A growing family of random graphs is called robust if it retains a giant component after percolation with arbitrary positive retention probability. We study robustness for graphs, in which new vertices are given a spatial position on the  \mbox{$d$-dimensional} torus  and are  connected to existing vertices with a probability favouring short spatial distances and high degrees. In this model of a
\mbox{scale-free} network with clustering  we can independently tune the power law exponent~$\tau$ of the degree distribution and
the rate $-\delta d$ at which the connection probability decreases with 
the distance of two vertices. We show that the network is robust if $\tau<2+\frac{1}{\delta}$, but fails to be robust if $\tau>3$.  In the case of one-dimensional 
space we also show that the network is not robust if  $\tau>2+\frac1{\delta-1}$. 
This implies that robustness of a scale-free network depends not only on its power-law exponent but also on its clustering features. 
Other than the classical models of scale-free networks our model is not locally tree-like, and hence  we need to develop novel methods for its study, 
including, for example, a surprising application of  the BK-inequality.  
\end{quote}

\vspace{0.3cm}

{\footnotesize
\noindent\emph{MSc Classification:}  Primary 05C80 Secondary 60C05, 90B15.

\noindent\emph{Keywords:} Spatial network, scale-free network, clustering, 
Barab\'asi-Albert model, preferential
attachment, geometric random graph, power law,  giant component,  robustness, phase transition, continuum percolation, 
disjoint occurrence, BK inequality.}

\vspace{-0.3cm}

\tableofcontents

\vspace{-1.5cm}

\section{Motivation}

Scientific, technological or social systems can often be described as complex networks of interacting components. Many of these networks 
have been empirically found to have strikingly similar topologies, shared features being that they are \emph{scale-free}, i.e.\ the degree distribution follows a power law, 
\emph{small worlds}, i.e.\ the typical distance of nodes is logarithmic or doubly logarithmic in the network size, or \emph{robust}, i.e.\ 
the network topology is qualitatively unchanged if an arbitrarily large proportion of nodes is removed from the network. Barab\'asi and Albert~\cite{Barabasi_Albert} 
therefore concluded fifteen years ago `that the development of large networks is governed by robust self-organizing phenomena that go beyond the particulars 
of the individual systems.'  They suggested a model of  a growing family of graphs, in which new vertices are added successively and connected to 
vertices in the existing graph with a probability proportional to their degree, and a few years later these features were rigorously verified in the work of
Bollob\'as and Riordan, see \cite{BRST01, BR03, BR04}.  
\medskip

In the years since the publication of~\cite{Barabasi_Albert} there have been many refinements of the idea of preferential attachment, introducing for 
example tunable power law exponents~\cite{Redner, DEHH09, DM09}, node fitness~\cite{Bianconi, Borgs, D10, DO14},  or spatial positioning of nodes~\cite{Flaxman, Aiello, Jordan}. Some of these refinements attempt to introduce or explain \emph{clustering}, the formation of clusters of nodes with an edge density significantly higher than in the overall network.
The phenomenon of clustering is present in real world networks but notably absent from most mathematical models of scale-free networks. The present paper investigates a \emph{spatial} network model, introduced in~\cite{JM14}, defined as a growing family of graphs in which a new vertex  gets a randomly allocated spatial position representing its individual features. This vertex then connects to every vertex in the existing graph independently,  with a probability which is a decreasing function of the spatial distance  of the vertices, the time, and  the inverse of the degree of the vertex. 
The relevance of this \emph{spatial preferential attachment model} lies in the fact that, while it is still a scale-free network governed by a simple rule of self-organisation, it  has been shown to exhibit  clustering.   The present paper investigates the problem of robustness and is probably the first rigorous attempt to understand the global topological structure of a self-organised scale-free network model with clustering.%
\medskip%

In mathematical terms, we call a growing family of graphs \emph{robust} if the critical parameter for vertex percolation is zero, which means that whenever vertices are deleted independently at random from the graph with a positive retention probability, a connected component comprising an asymptotically positive proportion of vertices remains.
For several scale-free models, including non-spatial preferential attachment networks, it has been shown  that robustness holds if  the power law exponent satisfies $\tau< 3$, see for example~\cite{BR03, DM13}. At a first glance one would maybe expect this behaviour to persist in the spatial model. It is known that robustness in scale-free networks relies on the presence of a hierarchically organised  core of vertices  with extremely high degrees, such that every vertex is connected to the next higher layer by a small number of edges, see for example~\cite{NR08}. Our analysis of the spatial model shows that, although a hierarchical core still exists if $\tau<3$, whether vertices in the core are sufficiently close in the graph distance to the next higher layer depends critically on the speed at which the connection probability decreases with spatial distance, and hence depending on this speed robustness may hold or fail.  The observation that robustness depends not only on the power law exponent but also on the clustering of a network appears to be new, though similar observations have been made for a long-range percolation model~\cite{DHH13}. Unlocking this phenomenon is the main achievement of this paper.
\medskip

The main structural difference between the spatial and classical model of preferential attachment is that the former exhibits \emph{clustering}. Mathematically this is measured in terms of a positive clustering coefficient, meaning that, starting from a randomly chosen vertex, and following two different edges, the probability that the two end vertices of these edges
are connected remains positive  as the graph size is growing. This implies in particular that local neighbourhoods of typical vertices in the 
spatial network do not look like trees. However, the main ingredient in almost every mathematical analysis of scale-free networks so far has 
been the approximation of these neighbourhoods by suitable random trees, see~\cite{BR09, DM11, BBCS14, EM14}.  As a result, the analysis of  spatial preferential attachment models requires a range of entirely new methods, which allow  to study the robustness of networks  without relying on the local tree structure that turned out to be so useful in the past. Providing these new methods is the main 
technical innovation in the present~work.
\pagebreak[3]

\section{The spatial preferential attachment model}

While spatial preferential attachment models may be defined in a variety of metric spaces,
we focus on homogeneous space represented by a $d$-dimensional torus of unit volume, given as
$\T_{1} =(-1/2,1/2]^{d}$ with the metric $\tmtextsf{d}_{1}$ given by
\[ \tmtextsf{d}_{1} (x,y) = \min   \big\{ \tmtextsf{d} (x,y+u):u \in \{
   -1,0,1\}^{d} \big\}, \mbox{ for } x,y \in \T_{1},\]
where $\tmtextsf{d}$ is the Euclidean distance on $\R^{d}$. Let $\X$ denote a homogeneous Poisson point process 
of finite intensity~$\lambda>0$ on $\T_{1} \times (0, \infty )$. A point $\x = (x,s)$ in $\X$ is a vertex $\x$,
born at time~$s$ and placed at position~$x$. 
Observe that, almost surely, two points of $\X$ neither have the same birth time nor the same position. 
We say that $(x,s)$ is {\emph{older}} than $(y,t)$ if $s<t$. For $t>0$, write $\X_{t}$ for
$\X \cap ( \T_{1} \times (0,t])$, the set of vertices already born at time~$t$. 
\medskip

We construct a growing sequence of graphs $(G_{t} )_{t>0}$, starting from
the empty graph, and adding successively the vertices in $\X$ when they are
born, so that the vertex set of $G_{t}$ equals~$\X_{t}$. Given the graph $G_{t-}$ at the time of birth 
of a vertex $\y = (y,t)$, we connect $\y$, independently of everything else, to each vertex $\x =(x,s)\in G_{t-}$,  
with probability
    \begin{equation}
      \label{construction} \varphi \left( \frac{t}{f( \indeg ( \x ,t-))}
      \tmtextsf{d}_{1} (x,y)^{d} \right) ,
    \end{equation}
where  $\indeg ( \x ,t-)$ is the \emph{indegree} of vertex $\x$, defined as the total number of edges between~$\x$ and younger vertices, at time~$t-$.  The model parameters in \eqref{construction} are
the {\emph{attachment rule} $f\colon \N \cup \{0\} \to (0, \infty )$, which is a nondecreasing function
regulating the strength of the preferential attachment, and the {\emph{profile
  function} $\varphi \colon [0, \infty ) \to (0,1)$, which is an integrable	nonincreasing function regulating the decay of
the connection probability in terms of the interpoint distance.%
\medskip%

The connection probabilities in~\eqref{construction} may look arcane at a first glance, but are in fact completely natural. 
To ensure that the probability of a new vertex connecting to its nearest neighbour does not degenerate, as $t\uparrow\infty$, it 
is necessary to scale $\tmtextsf{d}_{1} (x,y)$ by $t^{-1/d}$, which is the order of the distance of a point to its nearest neighbour 
at time~$t$. The linear dependence of the argument of $\varphi$ on time ensures that the expected number of edges connecting a new 
vertex to vertices of bounded degree remains bounded from zero and infinity, as $t\uparrow\infty$, 
as long as $\varphi$ is integrable on $[0, \infty )$, or equivalently $x \mapsto \varphi(\tmtextsf{d} (x,0)^d)$ is integrable on $\R^d$.
\medskip

The model parameters $\lambda$, $f$ and $\varphi$ are not independent. Indeed, if $\int \varphi(\tmtextsf{d}(x,0)^d) \, dx=\mu>0$, 
we can modify $\varphi$ to $\phi\circ (\mu \,\mathrm{Id})$ and $f$ to $\mu f$, so that the connection probabilities remain unchanged and
 \begin{equation}
    \label{normalisation} \int\varphi ( \tmtextsf{d} (x,0)^{d} ) \, dx=1.
  \end{equation}
Similarly, if the intensity of the Poisson point process $\X$ is $\lambda>0$, we can replace $\X$ by $\{(x,\lambda s) \colon (x,s)\in \X\}$ and 
$f$ by $\lambda f$, so that again the connection probabilities are unchanged and we get a Poisson point process of unit intensity. From now on
we will assume that both of these normalisation conventions are in place. 
\medskip

Under these assumptions the regime for the attachment rule~$f$ which leads to power law degree distributions is characterised by asymptotic linearity, i.e. 
$$\lim_{k\uparrow\infty}  \frac{f(k)}{k} = \gamma,$$
for some $\gamma>0$. We henceforth assume asymptotic linearity with the additional constraint that $\gamma<1$, which 
excludes degenerate cases with infinite mean degrees.
\pagebreak[3]

We finally assume that the profile function $\varphi$ is either regularly varying at infinity with index $-\delta$, for some $\delta >1$, or $\varphi$ decays 
quicker than  any regularly varying function.
 In the latter case we set $\delta = \infty$.  Intuitively, the bigger $\delta$, the stronger the clustering in the network. Our assumptions, in particular the assumption that 
$\varphi$ does not take the values 0 or 1, help us avoid some geometric constraints that are not of major interest. See Figure~1 for simulations of the spatial preferential attachment network indicative of the parameter dependence.%
\medskip%

\begin{figure}[ht]
\begin{center}
\includegraphics[scale=0.123]{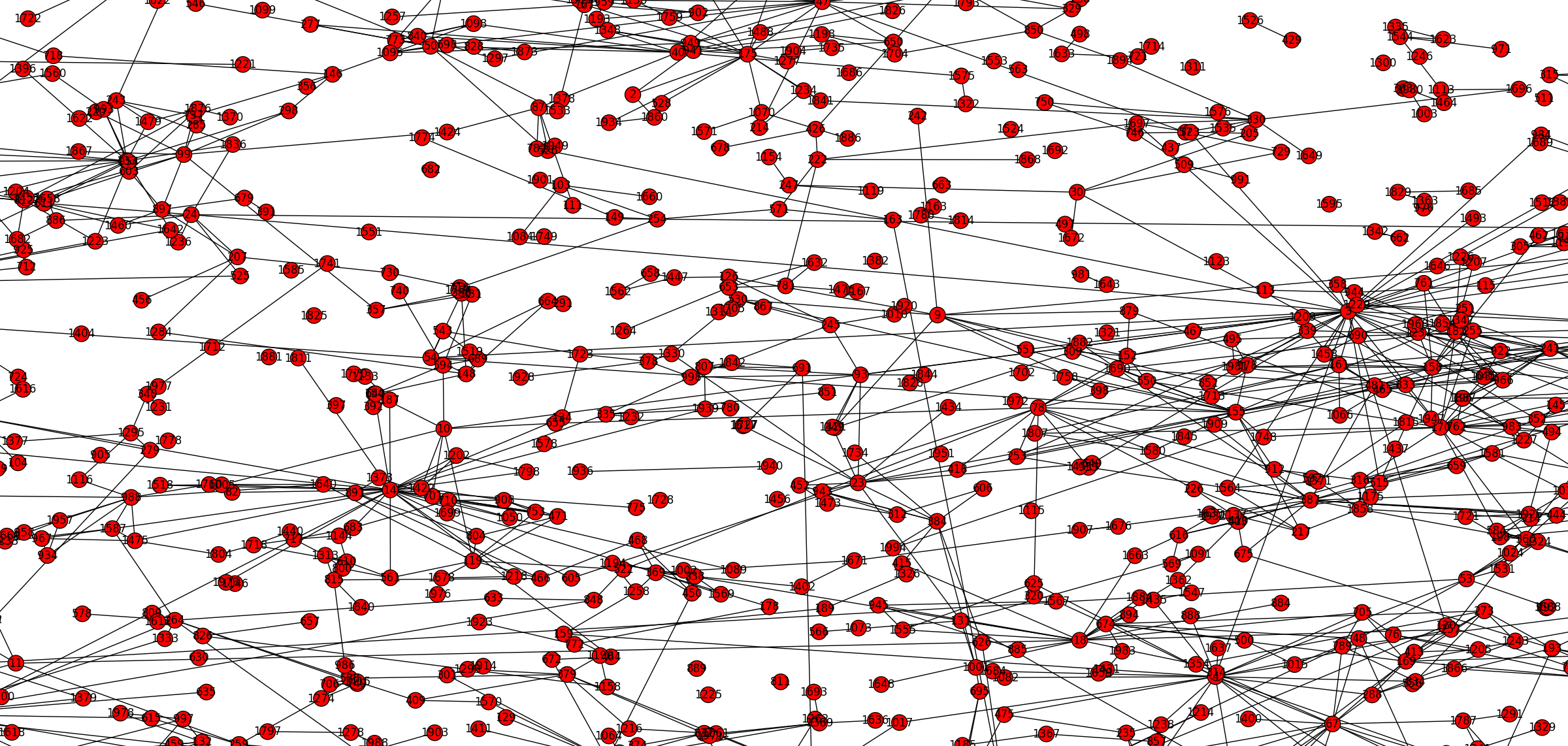}
\hspace{1mm}
\includegraphics[scale=0.123]{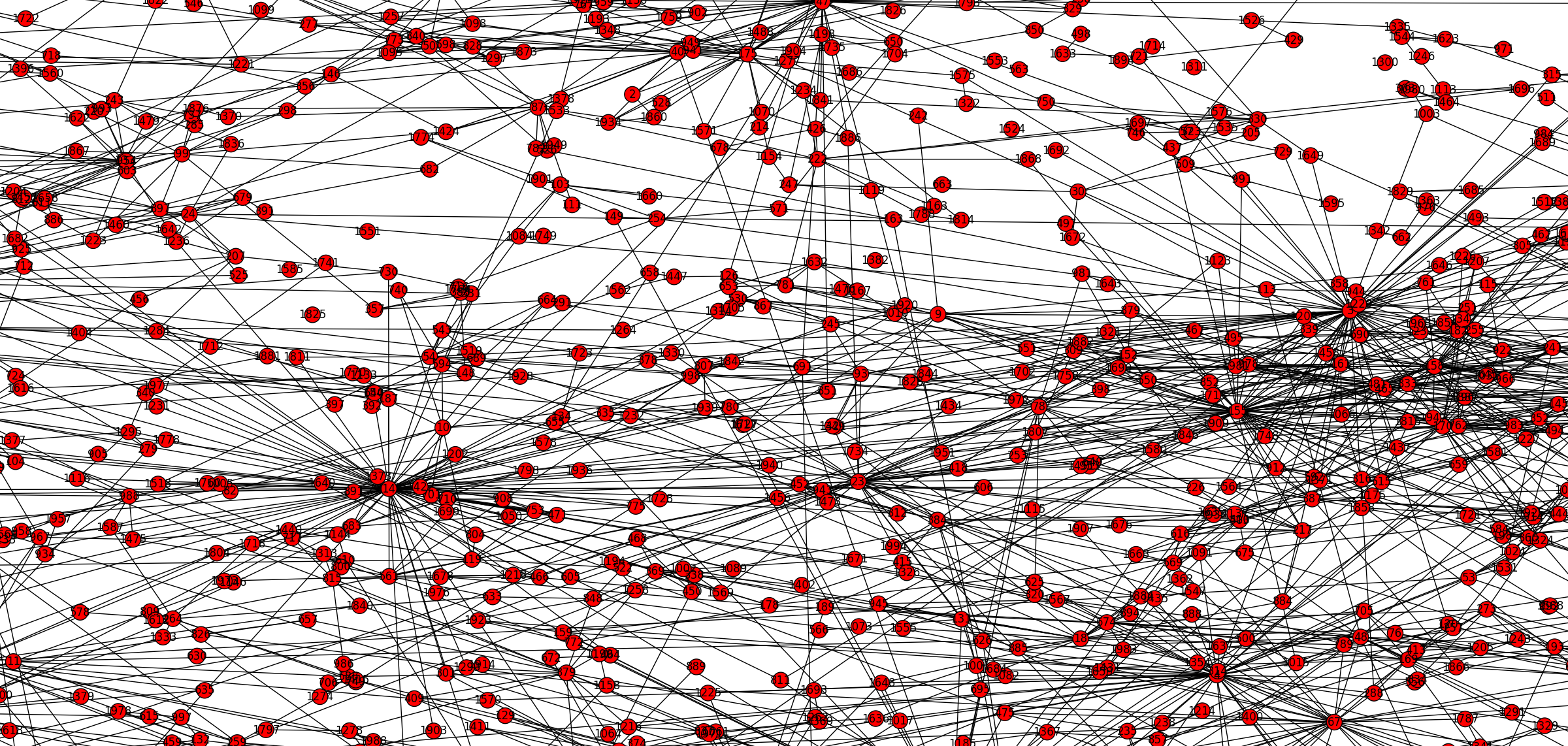}
\\[3mm]
\includegraphics[scale=0.123]{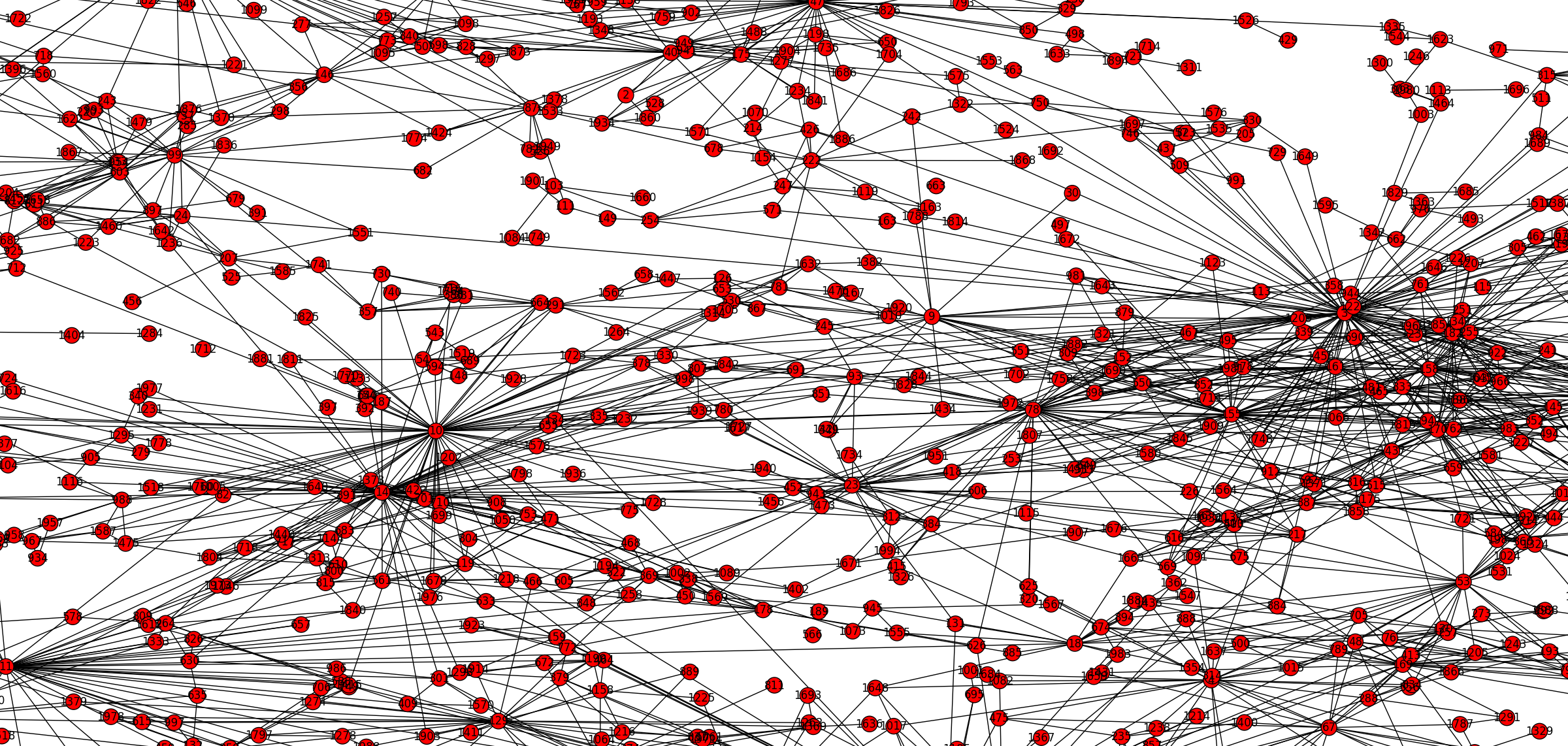}
\hspace{1mm}
\includegraphics[scale=0.123]{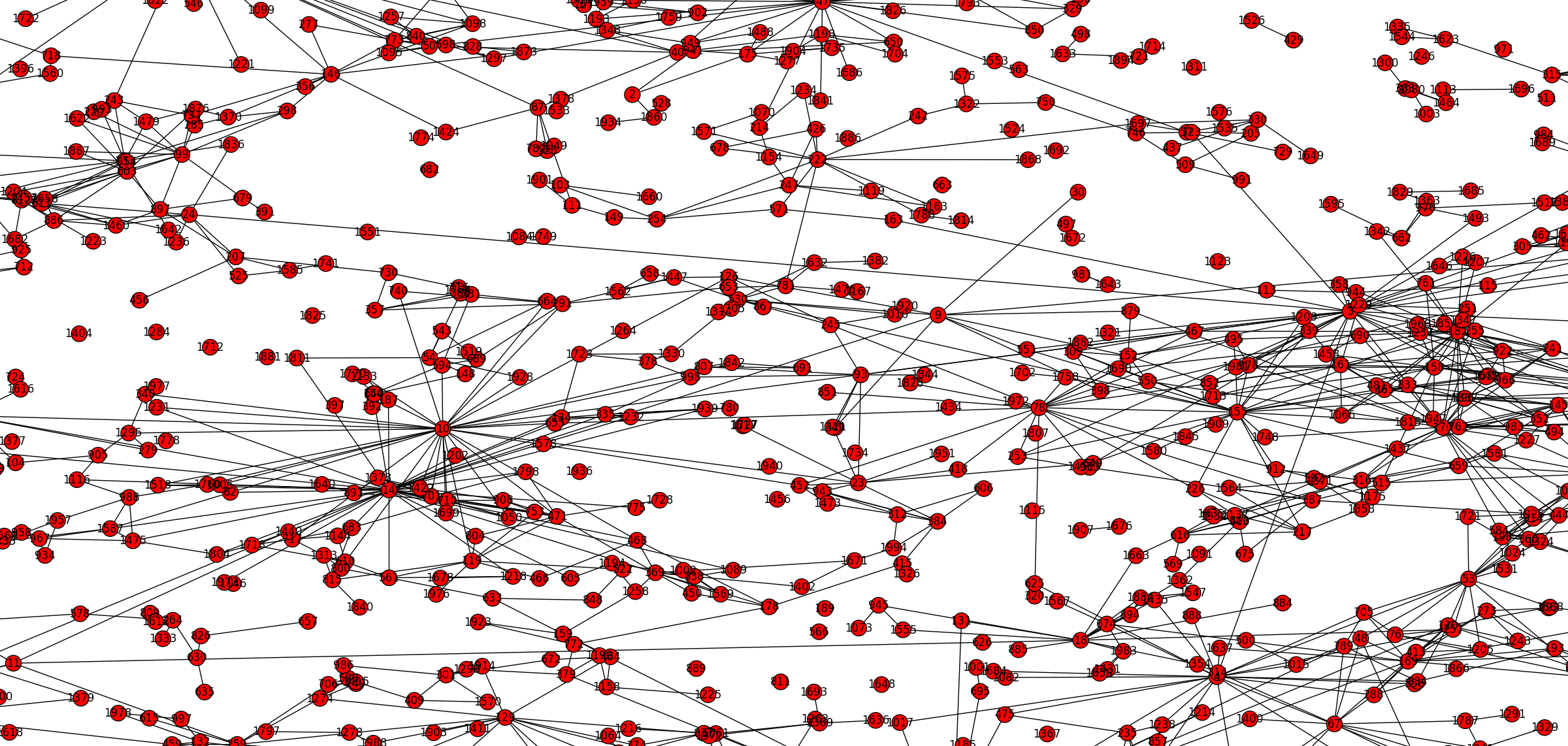}
\end{center}
\caption{Four simulations of the network in the case $d=2$, based on the same realisation of the Poisson process, with parameters (clockwise from top left) (a) $\gamma=0.5$ and $\delta=2.5$,
(b) $\gamma=0.75$ and $\delta=2.5$, (c)  $\gamma=0.5$ and $\delta=5$, (d) $\gamma=0.75$ and $\delta=5$. The pictures zoom into a 
typical part of the torus.}
\end{figure}

A similar spatial preferential attachment model was introduced in \cite{Aiello} and studied further in \cite{ JPW, Cooper}. In this model it is assumed that the profile functions has bounded support, more precisely $\phi =p \one_{[0,r]}$, for $p\in(0,1]$ and~$r$ satisfying~\eqref{normalisation}. This choice of profile function, roughly corresponding to the boundary case $\delta\uparrow\infty$, is too restrictive for the problems we study in this paper, as it turns out that robustness does not hold for any value of~$\tau$. There are also spatial long-range
percolation models which have a qualitatively similar behaviour to our networks, see for example~\cite{DHH13, DW13, Hirsch}, but these models are easier to analyse and the  methods 
of this analysis are quite different.%
\medskip%

Local properties of the spatial preferential attachment model were studied in~\cite{JM14}, where this model was first introduced. It is shown there, among other things, that
\begin{itemize}
\item The \emph{empirical degree distribution} of $G_{t}$ converges in probability to a deterministic limit~$\mu$. The probability measure~$\mu$ 
on~$\{0\}\cup\N$ satisfies
$$\mu(k) = k^{-(1+ \frac 1 \gamma)+ o(1)} \qquad \text{ as }\ \ k\uparrow \infty.$$
In other words, the network $(G_{t})_{t>0}$ is \emph{scale-free} with power-law exponent $\tau=1+ \frac 1 \gamma$, which can be tuned to take any value $\tau>2$. See~\cite[Theorem 1 and 2]{JM14}.\\[-1mm] \pagebreak[3]

\item The average over all vertices~$v\in G_t$ of the empirical local clustering coefficient at~$v$, defined as the 
proportion of pairs of neighbours of~$v$ which are themselves connected by an edge in $G_t$, converges in probability to a positive constant $c^{\rm av}_\infty>0$, called the
\emph{average clustering coefficient}. In other words the network $(G_{t})_{t>0}$ exhibits \emph{clustering}. See~\cite[Theorem 3]{JM14}.
\end{itemize}

\pagebreak[3]

\section{Statement of the main results}

We now address the problem of robustness of the network $(G_t)_{t>0}$ under percolation. Recall that the number of vertices of the graphs $G_t$, $t>0$, 
form a  Poisson process of unit intensity, and is therefore almost surely equivalent to $t$ as $t\uparrow\infty$. Let $C_t\subset G_t$ be the largest connected component in $G_t$  and denote by $|C_t|$ its size.  We say that the network has a \emph{giant component} if $C_t$ is of linear size or, more precisely, if
$$\lim_{\eps \downarrow 0} \limsup_{t\to \infty} \P\left(\frac{\vert C_t \vert} t\le \eps\right)=0.$$
We say it has \emph{no giant component} if $C_t$ has sublinear size or, more precisely, if
$$\liminf_{t\to \infty} \P\left(\frac{\vert C_t \vert} t\le \eps\right)=1 \mbox{ for any $\eps>0$.}$$
If $G$ is a graph with vertex set $\X$, and $p\in(0,1)$, we write 
$\tensor[^p]{G}{}$  for the random subgraph of $G$ obtained by Bernoulli
percolation with retention parameter $p$ on the vertices of $G$. We also use $\tensor[^p]{\X}{}$ for set of vertices surviving percolation.
The network $(G_t)_{t>0}$ is said to be \emph{robust} if, for any fixed $p\in(0,1]$, the network $(\tensor[^p]{G}{}_t)_{t>0}$ has a giant component
and \emph{non-robust} if there exists $p\in(0,1]$ so that $(\tensor[^p]G{}_t)_{t>0}$ has no giant component. 
\medskip

Our main result concern phases of robustness or non-robustness for the spatial preferential attachment network. In classical non-spatial preferential attachment models, there is a phase transition for robustness when $\tau$ crosses the critical value $3$, see~\cite{DM13}. It is easy to believe that the spatial structure does not help robustness. Our main result shows that in the spatial model robustness is still possible, but at least in the case $d=1$ the phase transition does not occur at $\tau =3$, but at a smaller value depending on $\delta$.
\medskip


\begin{theorem}
  \label{main} The spatial preferential attachment network $(G_t)_{t>0}$ is
  \begin{itemize}
    \item[(a)] \emph{robust} if $\gamma > \frac{\delta}{1+ \delta}$ or, equivalently, if $\tau <2+ \frac{1}{\delta}$;
    \item[(b)] \emph{non-robust} if $\gamma < \frac{1}{2}$ or, equivalently, if $\tau >3$.
  \end{itemize} 
  In the case of one space dimension, $d=1$, the network is
  \begin{itemize}   
    \item[(c)] \emph{non-robust} if $\gamma < \frac{\delta -1}{\delta}$ or, equivalently, if $\tau >2+ \frac{1}{\delta -1}$.
  \end{itemize}
\end{theorem}

\begin{figure}[ht]
\begin{center}
\includegraphics[scale=0.5]{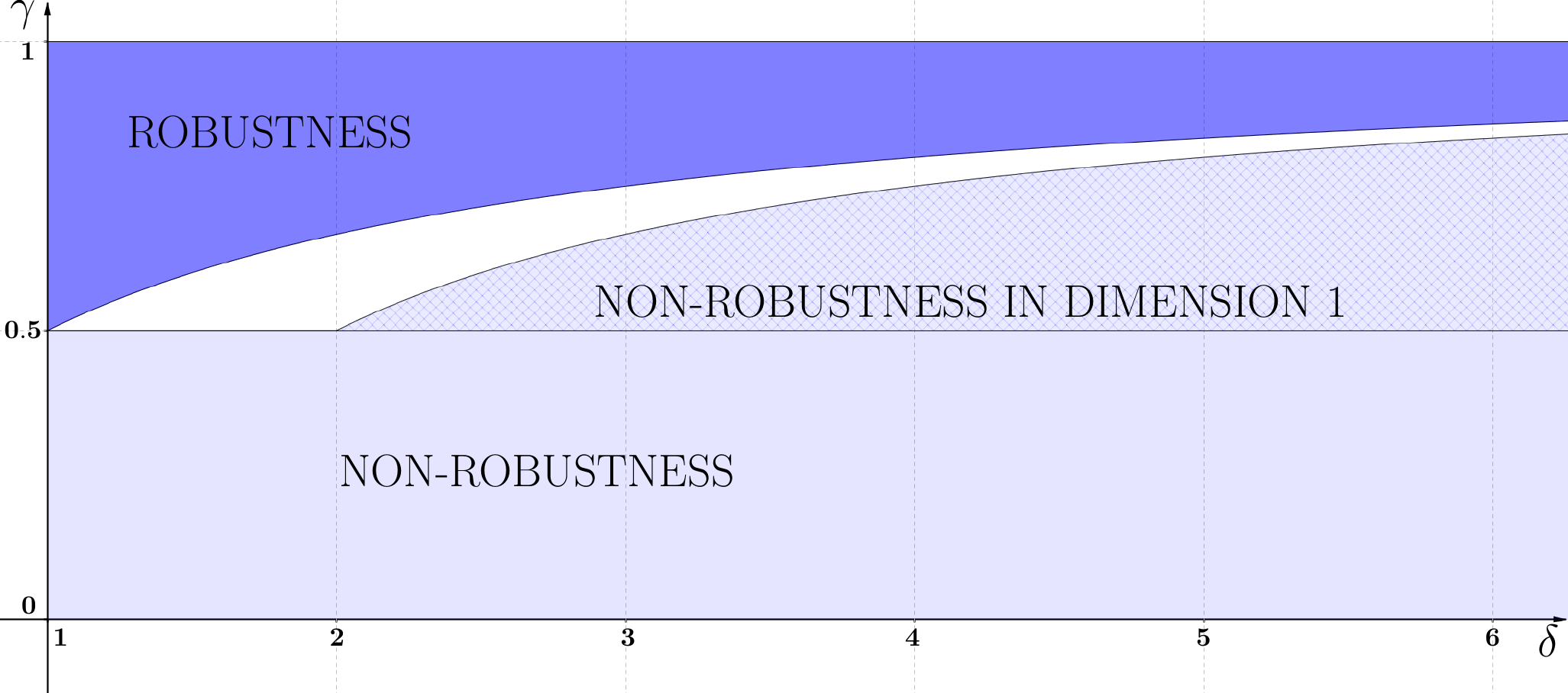}
\end{center}
\caption{The densely shaded area corresponds to the known robustness phase $\gamma>\delta/(\delta+1)$. The lightly shaded area corresponds to 
the known non-robustness phase. In between, no result of either robustness or non-robustness is known, though non-robustness is conjectured.
}
\end{figure}

\begin{rem}
For a suitable range of parameters, this seems to be the first instance of a  
scale-free network model which combines robustness with clustering features. 
We conjecture that nonrobustness occurs in any dimension if {\smash{$\gamma < \frac{\delta}{1+
\delta}$}}, and thus the critical value for $\tau$ equals {\smash{$2+
\frac{1}{\delta}$}}, but our proof techniques do not allow to prove this, see Figure~2 for a phase diagram.
\end{rem}

\pagebreak[3]

\begin{rem}
Our approach also provides heuristics indicating that in the robust phase $\delta(\tau-2)<1$ the typical distances in the robust giant component are asymtotically 
$$ (4+o(1)) \frac {\log \log t} {-\log(\delta (\tau-2))},$$
namely doubly logarithmic, just as in some nonspatial preferrential attachment models. The constant coincides with that of the nonspatial models in the limiting case 
$\delta\downarrow 1$, see~\cite{Dommers,DMM12}, and goes to infinity as $\delta(\tau-2)\to 1$. It is an interesting open problem to confirm these heuristics rigorously.
\end{rem}



\section{The limit model and proof strategies}

Before describing the strategies of our proofs, we briefly summarise the techniques developed in~\cite{JM14} in order to
describe the local neighbourhoods of typical vertices by a limit model. We will heavily rely on these techniques in the present paper.

\subsubsection*{Canonical representation}

We first describe a canonical representation of our network $(G_t)_{t>0}$. To this end, let $\X$ be a Poisson process of unit intensity 
on $\T_{1} \times (0, \infty )$, and endow the point process $\X\times \X$ with independent marks which are uniformly distributed 
on $[0,1]$. We denote these marks by $\V_{\x , \y}$ or $\V ( \x , \y )$, for $\x, \y\in\X$.%
\medskip

If $\Y\subset \T_{1} \times (0, \infty )$ is a finite set and $\W\colon\Y\times\Y\to[0,1]$ a map, 
we define a graph $G^{1}( \Y , \W )$ with vertex set $\Y$ by establishing edges in order of age of the younger endvertex.
An edge between $\x=(x,t)$ and $\y=(y,s)$, $t<s$, is present if and only~if
\begin{equation}
  \label{attachment_rule} \W ( \x , \y ) \le \varphi \left( \frac{s\,
  \tmtextsf{d}_{1} ( x , y )^{d}}{f(Z(\x,s-))} \right),
\end{equation}
where $Z(\x,s-)$ is the indegree of $\x$ at time $s-$. A realization of $\X$ and $\V$ then gives rise to
the family of graphs $(G_t)_{t>0}$ with vertex sets $\X_{t}=\X \cap ( \T_{1} \times (0,t])$, 
given by $G_t=G^{1} ( \X_{t} , \V )$,  which has the distribution of the spatial preferential attachment 
network.

\subsubsection*{Space-time rescaling}

The construction above can be generalised in a straightforward manner from $\T_{1}$ to the torus of volume $t$, 
namely $\T_{t}=(- \frac12\, t^{1/d}, \frac12\, t^{1/d}]^{d}$, equipped with its canonical torus metric~$\tmtextsf{d}_{t}$. The resulting
functional, mapping a finite subset $\Y\subset\T_{t} \times (0, \infty )$ and a map 
from $\Y \times\Y\to[0,1]$ onto a graph, is now denoted by $G^t$. 
\pagebreak[3]

We introduce the \emph{rescaling mapping}
\[ \begin{array}{llll}
     h_{t} : & \T_{1} \times (0,t] & \to & \T_{t} \times (0,1] ,\\
     & (x,s) & \mapsto & (t^{1/d} x,s/t)
   \end{array} \]
which expands the space by a factor $t^{1/d}$, the time by a factor $1/t$. The
mapping $h_{t}$ operates on the set $\X$, but also on $\V$, by
$h_{t} ( \V )_{h_{t} ( \x ) ,h_{t} ( \y )} \assign \V_{\x , \y}$.
The operation of $h_{t}$ preserves the rule~{\eqref{attachment_rule}},
and it is therefore simple to verify that we have
\[ G^{t} (h_{t} ( \X_{t} ),h_{t} ( \V )) =h_{t} (G^{1} ( \X_{t} , \V )) =h_{t}
   (G_{t} ) , \]
that is, it is the same to construct the graph and then rescale the picture,
or to first rescale the picture, then construct the graph on this rescaled
picture. Observe also that $h_{t} ( \X_{t} )$ is a {\PPP} of intensity $1$ on
$\T_{t} \times (0,1]$, while $h_{t} ( \V )$ are independent marks attached to 
the points of $h_{t} ( \X_{t} )\times h_{t} ( \X_{t} )$ which are
uniformly distributed on $[0,1]$.

\subsubsection*{Convergence to the limit model}

We now denote by $\X$ a {\PPP} with unit intensity on $\R^d \times (0,1]$, and endow
the points of $\X\times\X$ with independent marks $\V$, which are uniformly distributed on~$[0,1]$. 
For each $t>0$, identify $(- \frac12\, t^{1/d}, \frac12\, t^{1/d} /2]^{d}$ and
$\T_{t}$, and write $\X^{t}$ for the restriction of $\X$ to $\T_{t} \times
(0,1]$, and $\V^{t}$ for the restriction of $\V$ to $\X^{t}\times\X^{t}$. 
In the following, we write $G^t$ or $G^{t} ( \X , \V )$ for $G^{t} (
\X^{t} , \V^{t} )$. We have seen that  for fixed $t \in (0, \infty )$, the graphs $G^t$ and 
$h_{t} (G_{t} )$ have the same law. Thus any results of robustness
we prove for the network $(G^t)_{t>0}$ also hold for the network $(G_t)_{t>0}$. 
It was shown in~\cite[Proposition~5]{JM14} that, almost surely, the graphs $G^t$ converge  to a 
locally finite graph $G^\infty= G^{\infty} (\X, \V)$, in the sense that the neighbours of any given vertex 
$\x\in\X$ coincide in $G^t$ and in $G^\infty$, if $t$ is large enough. 
It is important to note the fundamentally different behaviour of the processes
$(G^t)_{t>0}$ and $(G_t)_{t>0}$. While in the former the degree of any fixed vertex stabilizes, in the latter
the degree of any fixed vertex goes to $\infty$, as $t\uparrow\infty$. We will exploit the convergence of $G^t$ to $G^\infty$ 
in order to  decide the robustness of the finite graphs $G^t$, and ultimately $G_t$, 
from properties of the limit model $G^\infty$.

\subsubsection*{Law of large numbers}

We now state a limit theorem for the graphs \smash{$\tensor[^p]{G}{^t}$} centred in a randomly chosen point. To this end we denote by $\tensor[^p]{\tmP}{}$ the law of $\X , \V$ together with independent Bernoulli percolation with retention parameter~$p$ on the points of~$\X$. For any $\x\in\R^d\times(0,1]$ we denote by $\tensor[^p]{\tmP}{}_{\x}$  the \emph{Palm measure}, i.e.\ the law $\tensor[^p]{\tmP}{}$ conditioned on the event  
$\{ \x\in\!\tensor[^p]{\X}{} \}$. Note that by elementary properties of the Poisson process this conditioning simply adds the point $\x$ to $\tensor[^p]{\X}{} $ and 
independent marks $\V_{\x, \y}$ and $\V_{\y, \x}$, for all $\y\in\X$, to $\V$.  We also write $\tensor[^p]{\mathbbm E}{} _{\x}$ for the expectation under $\tensor[^p]{\tmP}{}_{\x}$. 
Let $\xi = \xi \left(\x,G \right)$ be a bounded functional of a locally-finite graph $G$ with vertices in $\mathbbm{R}^d \times ( 0,1]$  and a vertex~$\x\in G$, which is invariant under translations of~$\R^d$.
Also, let \smash{$\xi_t = \xi_t \left(\x,G \right)$} be a bounded family of functionals of a graph $G$ with vertices in $\T_t \times ( 0,1 ]$ and a vertex $\x\in G$, invariant under translations of the torus. 
We assume that, for $U$ an independent uniform random variable on $(0,1]$, we have 
that $\xi_t((0,U), \tensor[^p]{G}{^t})$ converges to $\xi((0,U), \tensor[^p]{G}{^\infty})$ in $\tensor[^p]{\P}{}_{(0,U)}$-probability.
Then, in $\tensor[^p]{\P}{}$-probability,
\begin{equation}\label{LLN}
  \frac{1}{t} \sum_{\x \in \tensor[^p]{\X}{^t}} \xi_t \left( \x
  , \tensor[^p]{G}{^{t}} \right) \underset{t \rightarrow \infty}{\longrightarrow}
 p\, \int_0^1\tensor[^p]{\E}{}_{(0,u)}[\xi((0,u), \tensor[^p]{G}{^\infty})] \, du.
\end{equation}
This law of large numbers is a minor modification of the one given in~\cite[Theorem~7]{JM14}, which covers the case $p=1$.
\bigskip

\pagebreak

\subsection{Robustness: strategy of proof}\label{rob} 

\subsubsection*{Existence of an infinite component in the limit model} We first show that, under the assumptions that 
\smash{$\gamma>\frac{\delta}{1+\delta}$} and $p\in(0,1]$,  the percolated limit model \smash{$\tensor[^p]{G}{^\infty}$} has an infinite connected component. This uses the 
established strategy of the {hierarchical core}. The young vertices,  born after time $\frac12$, are called connectors. 
Fix $\alpha\in(1,\frac{\gamma}{\delta(1-\gamma)})$. Starting from a sufficiently old vertex~\smash{$\x_0\in\tensor[^p]{G}{^\infty}$}, we establish an infinite chain \smash{$(\x_k)_{k\ge 1}$} of vertices $\x_k=(x_k, s_k)$ such that $s_k<s_{k-1}^\alpha$, i.e. we move to increasingly older vertices,
and $\x_{k-1}$ and $\x_{k}$ are connected by a path of length two, using a connector as a stepping stone.

\subsubsection*{Transfer to finite graphs using the law of large numbers} To infer robustness of the network \smash{$(G^t)_{t>0}$}  from the 
behaviour of the limit model we use~\eqref{LLN} on the functional $\xi_t(\x, G)$ defined as the indicator of the event that there is a path in~$G$ connecting $\x$ to the oldest vertex of~$G$. We denote by $\xi(\x, G)$ the indicator of the event that the connected component of $\x$ is infinite and let
\begin{equation}\label{thetadef}
\tensor[^p]{\theta}{} := \int_0^1\tensor[^p]{\P}{}_{(0,u)}\big\{ \mbox{ the component of } (0,u) \mbox{ in }
 \tensor[^p]{G}{^\infty} \mbox{ is infinite }\big\} \, du.
\end{equation} 
If $\lim \xi_t((0,U), \tensor[^p]{G}{^t})=\xi((0,U),\tensor[^p]{G}{^\infty})$ in probability, then the law of large numbers~\eqref{LLN}~implies
$$ \frac 1 {t} \sum_{\x \in \tensor[^p]{\X}{^t}} \xi_t(\x, \tensor[^p]{G}{^t}) \longrightarrow p\, \tensor[^p]{\theta}{}.$$
The sum on the left is the number of vertices in $\tensor[^p]{G}{^t}$ connected to the oldest vertex, and we infer that this number grows linearly in $t$ 
so that a giant component exists in $(\tensor[^p]{G}{^t})_{t>0}$. This implies that $({G}{^t})_{t>0}$ and hence $({G}{_t})_{t>0}$ is a robust network. 
%
However, while it is easy to see that
$\limsup_{t\uparrow\infty} \xi_t((0,U), \tensor[^p]{G}{^t}) \leq \xi((0,U),\tensor[^p]{G}{^\infty}),$
checking that
\begin{equation}\label{hardpart}
\liminf_{t\uparrow\infty} \xi_t((0,U), \tensor[^p]{G}{^t}) \geq \xi((0,U),\tensor[^p]{G}{^\infty}),
\end{equation}
is the difficult part of the argument.

\subsubsection*{The geometric argument} The proof of~\eqref{hardpart} is the most technical part of the paper. We first look at the finite graph~$\tensor[^p]{G}{^t}$ and establish the existence of a {core} of old and well-connected vertices, which includes the oldest vertex. Any pair of vertices in the core are connected by a path with a bounded number of edges, in particular all vertices of the core are in the same connected component.
This part of the argument is similar to the construction 
in the limit model. We then use a simple continuity argument to establish that if the vertex~$(0,U)$ is in an infinite component in the limit model, then it is also in an infinite component for the limit model based on a Poisson process $\X$ with a slightly reduced intensity.  In the main step we show that under this assumption the vertex $(0,U)$ is connected in $\tensor[^p]{G}{^t}$ with reduced intensity to a moderately old vertex. In this  step we have to rule out explicitly the possibilities that the infinite component of  $\tensor[^p]{G}{^\infty}$ either avoids the set of eligible moderately old vertices, or connects to them only by a path which moves very far away from the origin.  The latter argument requires good control over the length of edges in the component of $(0,U)$ in $\tensor[^p]{G}{^\infty}$. Once the main step is established, we can finally use the still unused  vertices, which form a Poisson process with small but positive intensity, to connect the moderately old vertex we have found to the core by means of a classical {sprinkling argument}.
\medskip

\pagebreak[3]

In Section~5 we will carry out this programme and prove robustness. In fact, we shall prove that under our hypothesis $\gamma>\frac{\delta}{1+\delta}$ a stronger statement holds, see Proposition~\ref{finerob}, which also implies that the size of the second largest connected component in $\tensor[^p]{G}{}_t$ does not grow linearly. In other words, we will see that in this regime the network has a \emph{unique} giant component.
\pagebreak[3]

\subsection{Non-robustness: strategy of proof}\label{nonrob}

\subsubsection*{Using the limit model} If $\gamma<\frac12$ it is very plausible that the spatial preferential attachment network is non-robust, as the classical models with the same power-law exponents are non-robust~ \cite{BR03, DM13} and it is difficult to see how the spatial structure could help robustness. We have not been able to use this argument for a proof, though, as our model cannot be easily dominated by a non-spatial model with the same power-law exponent. Instead we use a direct approach, which turns out to yield non-robustness also in some cases where $\gamma>\frac12$. The key is again the use of the limit model, and in particular the law of large numbers. We apply this now to the functionals $\xi^{\ssup k}(\x, G)$ defined as the indicator of the event that the connected component of~$\x$ has no more than~$k$ vertices. Clearly, $\lim_{t\uparrow\infty} \xi^{\ssup k}((0,U), \tensor[^p]{G}{^t})= \xi^{\ssup k}((0,U),\tensor[^p]{G}{^\infty}),$ and therefore
\begin{equation} \label{finitecomponents}
\frac 1 {pt} \sum_{\x \in \tensor[^p]{\X}{^t}} \xi^{\ssup k}(\x, \tensor[^p]{G}{^t}) \longrightarrow \int_0^1\tensor[^p]{\E}{}_{(0,u)}[\xi^{\ssup k}((0,u), \tensor[^p]{G}{^\infty})] \, du.
\end{equation}
The left hand side is asymptotically equal to the proportion of vertices in $\tensor[^p]{G}{^t}$ which are in components no bigger than~$k$. As $k\to \infty$
the right hand side converges to $1-\tensor[^p]{\theta}{}$. Hence if $\tensor[^p]{\theta}{}=0$ for some $p>0$, then 
$({G}{^t})_{t>0}$ and hence $({G}{_t})_{t>0}$ is non-robust. It is therefore sufficient to show that, for some sufficiently small~$p>0$, 
there is no infinite component in the percolated limit model~$\tensor[^p]{G}{^\infty}$.

\subsubsection*{Positive correlation between edges}
We first explain why a na\"{i}ve first moment calculation fails. If $(0,U)$ has positive probability of belonging to an infinite component of $\tensor[^p]{G}{^\infty}$ then, with positive probability, we could find an infinite self-avoiding path in  $\tensor[^p]{G}{^\infty}$ starting from $\x_0=(0,U)$.
A direct first moment calculation would require to give a bound on the probability of the event
$\{\x_0{\leftrightarrow}\x_1{\leftrightarrow }\cdots{\leftrightarrow}\x_n\}$ that a sequence $(\x_0,\ldots, \x_n)$ of distinct points~$\x_i=(x_i, s_i)$ conditioned to be in $\X$ forms a path in ${G}{^\infty}$. If this estimate allows us to bound the expected number of paths of length~$n$ in ${G}{^\infty}$ starting in $\x_0=(0,U)$ by $C^n$, for some constant $C$, we can infer with Borel-Cantelli that, if $p<1/C$, almost surely there is no 
arbitrarily long self-avoiding paths in~$\tensor[^p]{G}{^\infty}$. For a variety of non-spatial models the event  $\{\x_0{\leftrightarrow}\x_1{\leftrightarrow }\cdots{\leftrightarrow}\x_n\}$ can be decomposed into independent, or negatively correlated, events of the form $\{\x_i {\leftrightarrow} \x_{i+1}\}$, or $\{\x_{i-1} {\leftrightarrow} \x_i {\leftrightarrow} \x_{i+1}\}$ with $s_i< s_{i-1}, s_i$, the probability of which can be easily estimated, see for example~\cite{DMM12}. For spatial networks however such a decomposition is not possible.
Indeed, the events $\{\x_j{\leftrightarrow}\x_{j+1}\}$ and $\{\x_k{\leftrightarrow}\x_{k+1}\}$ are not independent if the interval $I=(s_j,s_{j+1})\cap (s_k,s_{k+1})$ is nonempty, because the existence of a vertex 
in $\X\cap (\R^d\times I)$ which is relatively close to both $x_j$ and $x_k$ is likely to connect to both of these vertices and make their indegrees grow simultaneously. Observing that all the events $\{\x_k{\leftrightarrow}\x_{k+1}\}$ are increasing in the Poisson point process $\X$, we can argue by Harris' inequality that they are positively correlated. Because the positive correlations  play against us, it seems impossible to give an effective upper bound on the probability of a 
long sequence to be a path, therefore making this first moment calculation impossible.

\subsubsection*{Quick paths, disjoint occurrence, and the BK inequality}
As a solution to this problem we develop the concept of {quick paths}. Starting from a sequence 
$(\x_0, \ldots, \x_n)$ in $\tensor[^p]{G}{^\infty}$, with $\x_0=\x$ and $\x_n=\y$ we construct a new sequence  
$(\z_0, \ldots, \z_m)$, with  $\z_0=\x$ and $\z_m=\y$, such that at least half 
of the points  are in $\tensor[^p]{G}{^\infty}$, and the remaining ones are in ${G}{^\infty}$. This sequence, called a quick path, has the property that the event $\{\z_0{\leftrightarrow} \cdots{\leftrightarrow}\z_m\}$ can be split into smaller parts, in the sense that it implies
the \emph{disjoint occurrence} of events $\{\z_{i}{\leftrightarrow}\cdots{\leftrightarrow}\z_{i+4}\}$ involving five or fewer consecutive vertices of the sequence. The concept of disjoint occurrence is due to van den Berg and Kesten, and the famous BK-inequality states that the probability of events occurring disjointly is bounded by the probability of their product. It is tedious, but not hard, to estimate the probability of these paths of length no more than five, and the estimate produces the necessary bounds to complete the argument.%
\medskip%

\pagebreak[3]

Instead of defining quick paths and disjoint occurrence here, we just give a flavour by  
showing how to deal with a path $\x_0{\leftrightarrow}\x_1{\leftrightarrow}\x_2$, for 
$\x_0, \x_1, \x_2\in \tensor[^p]{G}{^\infty}$ and $s_0, s_2<s_1$. To move to the quick path we 
let $\z_0=\x_0$, $\z_2=\x_2$ and replace the vertex $\x_1\in \tensor[^p]{G}{^\infty}$ by the 
oldest vertex  $\z_1=(z_1, u_1) \in {G}{^\infty}$ such that $s_0, s_2<u_1$ and
$\z_0{\leftrightarrow}\z_1{\leftrightarrow}\z_2$. Now any vertex $\z_1'=(z_1', u_1') 
\in {G}{^\infty}$ with $s_0, s_2< u_1'<u_1$ can only influence the indegree  of either 
$\z_0$ or $\z_2$ at time $u_1$, but never both. This means, loosely speaking, that 
$\z_0{\leftrightarrow}\z_1{\leftrightarrow}\z_2$ being a quick path implies the disjoint occurrence of the events 
$\{\z_0{\leftrightarrow}\z_1\}$ and $\{\z_1{\leftrightarrow}\z_2\}$. %
\medskip%
\pagebreak[3]

In Section~6 we will carry out this programme and prove non-robustness if $\gamma<\frac12$, or if $d=1$ and 
\smash{$\gamma<\frac{\delta-1}{\delta}$.} Some auxiliary lemmas used in various parts 
of our proofs have been postponed to an appendix, see Section~7.

\subsubsection*{Summary of standard notation}

We use the Landau symbols \smash{$o (t)$, $O (t)$}, and $\Omega (t)$. If \smash{$(A (t))_{t>0}$} is a family of events, we say $A (t)$ holds with high probability, or {\whp}$(t)$, if the probability of $A (t)$ goes to 1, as $t\uparrow\infty$. We say $A (t)$ holds with extreme probability, or {\wep}$(t)$, if it
holds with probability at least $1-\exp(- \Omega (\log^{2} t) )$, as $t \uparrow
\infty$. When the parameter is clear, we write \whp or \wep for {\whp}$(t)$ or {\wep}$(t)$. Observe that, if $(A(t)_n)_{n\ge0}$ is 
a sequence of events that simultaneously hold {\wep}$(t)$, in the sense that $\sup_n \text{Prob}(A(t)_n^{\rm c})=\exp(- \Omega (\log^{2} t) )$, 
as $t \uparrow \infty$, and $n(t)$ satisfies $n(t)=o(t^\varkappa)$ for some $\varkappa>0$,  then $ \bigcap_{k\le n(t)} A(t)_{k}$ holds {\wep}$(t)$. Informally, the intersection of a polynomial number of events, which each holds {\wep}, holds~{\wep}.

\section{Proof of robustness}


In the following three subsections we study percolation on the infinite graph $G^\infty$. The giant component for the 
sequence $(G^t)_{t\geq 0}$ of finite graphs is studied in Subsection 5.4. 

\subsection{Infinite component of the infinite graph}\label{uniqueness_part}

In this section we prove that the infinite graphs $\tensor[^p]{G}{^\infty}$ cannot contain more than one infinite component.

\begin{prop}
 In the graphs
$\tensor[^p]{G}{^{\infty}},$ for $p\in(0,1]$, the number
  of infinite components is always either almost surely equal to zero, or almost surely
  equal to one.
\end{prop}


The analogue of this proposition for percolation on the integer lattice, or
the Poisson random connection model in $\R^{d}$, is known for some time. 
Our proof follows the classical technique of Burton and Keane, see~\cite{BuKe}. 
We focus on the case $p=1$, the other cases being similar.
The first step of this proof is to use the ergodicity of the model to deduce that the number of
infinite components is an almost sure constant in $\N \cup \{0, \infty\}$. The
second step is to ensure that this constant cannot be a finite number $k \ge
2$. Informally, if the graph contains $k \ge 2$ infinite components and if
$\x$ and $\y$ are two vertices belonging to 2 different infinite components,
then, sampling again the random variable $\V_{\x , \y}$ gives us a
positive probability of connecting these two components, hence decreasing the
number of infinite components by (at least) one, leading to a contradiction.
{\medskip}

We focus on the last step, which is to ensure that the number of connected components
cannot be infinity. We suppose, for contradiction, that
$$(H1)\hspace{5mm}\text{almost surely, the graph }G^\infty\text{ contains infinitely many infinite components.}$$
We say a vertex $\x \in G^{\infty}$ is a
\emph{trifurcation} if it is linked to at least three other vertices $\x_{1}, \x_{2}, \x_{3}\in G^{\infty}$, 
so that if $\x$ (and all its adjacent edges) is removed, then $\x_{1}$, $\x_{2}$, $\x_{3}$ are in three
different infinite components of $G^{\infty}$. Note that prior to the removal of $\x$ the vertices
$\x_{1}$, $\x_{2}$, $\x_{3}$ are all in the same
infinite component as they are all connected to~$\x$.

\begin{lemma}
  \label{trifurcations}If $(H1)$ is satisfied, then
$\displaystyle\int_0^1 du \, \P_{(0,u)}\big(\mbox{$(0,u)$ is a trifurcation}\big) >0.$
\end{lemma}

\begin{proof}
We write $\Pp(du \, d\omega)=\int_0^1 du\, \P_{(0,u)}(d\omega)$ for the underlying measure and note that $U=U(u,\omega)=u$  is the uniformly distributed birth time of the vertex located at the origin. Recall that $G^{\infty} =G^{\infty} ( \X , \V )$ and abbreviate
$G^{\infty}_{0} :=G^{\infty}  ( \X -\{(0,u)\},  \V )$, so that the law of $G^\infty_0$  under $\Pp$ is the same as the law 
of $G^\infty$ under $\P$. Observe that the conditional probability given $G_0^\infty$ and $U$, that the vertex located at the origin has degree 0 in $G^\infty$, is almost surely in $(0,1)$. This observation uses in particular the fact that $\phi$ does not take the values zero or one. Similarly, the conditional probability that the neighbouring vertices of $(0,U)$ in $G^\infty$ are exactly some given $\x_1$, $\x_2$, $\x_3$ in $G_0^\infty$, is also almost surely in $(0,1)$.
\smallskip

Assuming~$(H1)$, the graph $G^{\infty}_{0}$ contains almost surely infinitely many infinite components, and we may then specify arbitrarily three vertices $\x_1$, $\x_2$, $\x_3$ belonging to three different ones. Then, under $\Pp$, there is positive probability that $(0,U)$ is older than these three vertices, and connected by an edge to exactly these three vertices. If this happens, then the presence of the new vertex 
$(0,u)$ and the three edges linking to $\x_{1}$, $\x_{2}$ and $\x_{3}$ does not change the
indegree of any of the other vertices. In other words, it does not interfere with the rest of the graph. 
The graphs $G^{\infty}_{0}$ and $G^{\infty}$ are exactly the same except that the latter contains one vertex and three edges
more. Therefore $(0,u)$ is a trifurcation of~$G^{\infty}$. 
\end{proof}

By stationarity the expected number of trifurcations in the ball $\B (0,t)$, centered at the origin and of 
radius~$t$, is proportional to the volume of the ball, that is to \smash{$t^{d}$}. Let $E_{t}$ be the set of 
edges connecting a vertex $\x$ inside to a vertex $\y$ outside of $\B (0,t)$.

\begin{lemma}
  The cardinality of $E_{t}$ exceeds the number of trifurcations in $\B (0,t)$ by at least~2.
\end{lemma}

\begin{proof}
To each trifurcation $\x = (x,s)$ with $x \in \B (0,t)$, we can associate a
partition of $E_{t}$ into three sets such that two edges in $E_{t}$ that are
not in the same set are not in the same component of the graph with $\x$
and its incident edges removed. We obtain a compatible collection of
partitions in the sense of Burton-Keane and the result follows accordingly.
\end{proof}

By the lemma, the expected number of elements in~$E_{t}$ is at least proportional to $t^{d}$.
But the expected degree of a given vertex (with birth time
uniform in $(0,1]$) is finite, and so the expected number of neighbours at
distance at least $L$ is decreasing to zero, as $L \to \infty$. It follows that the 
expectation of $E_{t}$ must be $o (t^{d} )$, and
we obtain a contradiction.

\subsection{Continuity of the density of the infinite component}

In this subsection, we are interested in the continuity properties of $\tensor[^p]{\theta}{}$, as defined in~\eqref{thetadef}, 
with respect to the parameters of the model. To this end we now
suppose that $\tensor[^p]{\tmP}{}$ provides a consistent family of Poisson point processes $(_\lambda\X)_{0\le \lambda\le1}$ of intensity~$\lambda$ so that, for $\lambda<\mu$, we have $_\lambda\X\subset$ $\!_{\mu}\X$ and hence  $G^t(_\lambda\X, \V)$ is a subgraph of  $G^t(_{\mu}\X, \V)$.  For any $\x\in\R^d\times(0,1]$ we denote by $\tensor*[^p_\lambda]{\tmP}{}_{\x}$  the  law $\tensor[^p]{\tmP}{}$ conditioned on the event  
$\{ \x\in\!\tensor*[^p_\lambda]{\X}{} \}$. Denoting $_\lambda G^t:=G^t(_\lambda\X, \V)$ we again obtain convergence to a limit graph~$_\lambda G^\infty$. 
We let $\tensor*[_\lambda^p]{G}{^\infty}$ be the percolated limit graph. Recalling that $\xi ( \x ,G)$ is the indicator of the event that the connected component of~$\x$ in $G$ is infinite, we define
$$\tensor*[^p_\lambda]{\theta}{}:=
\int_0^1\tensor*[^p_\lambda]{\E}{}_{(0,u)} \big[ \xi\big( (0,u),  \tensor*[^p_\lambda]{G}{^\infty}\big) \big] \, du.$$

\begin{proposition}\label{con} For fixed $p\le1$, the function $\lambda\mapsto \tensor*[^p_\lambda]{\theta}{}$ is non-decreasing, right-continuous, and left-continuous everywhere except possibly at 
\smash{$\lambda_c(p):=\sup\{\lambda \colon\tensor*[^p_\lambda]{\theta}{}=0\}$.}
\end{proposition}

\begin{rem}
A similar result holds for the function $p\mapsto \tensor[^p]{\theta}{}$, and is a variant of well-known results of percolation theory.
\end{rem}

\begin{proof}[Proof of Proposition~\ref{con}.]
We only prove left-continuity at $\lambda>\lambda_c$ in the case $p=1$, as the other parts of the statement will not be used in the sequel.
Fix $\lambda_c<\mu<\lambda$, and let~$C$ be the infinite component of $_\mu G^{\infty}$. On the event that $(0,u)$ is in the infinite component of 
$_\lambda G^{\infty}$, it is connected in  this graph to a vertex of $C$. Thus there exists $k\in\N$ such that the
 $k$-neighbourhood of $(0,u)$ in the graph $_{\lambda}G^{\infty}$, i.e.\ the set of vertices of the graph within graph distance $k$ of $(0,u)$, 
intersects $C$.
Given~$k$, with probability one, the $k$-neighbourhood of $(0,u)$ is the same in $_{\lambda}G^{\infty}$ and in $_{\lambda'}G^{\infty}$ for $\lambda'$ close enough to $\lambda$. Hence $ \xi ( ( 0,u ) ,G^{\infty}(_{\lambda'}\X, \V) )$ converges almost surely, and hence in probability, to  $\xi ( ( 0,u ) ,G^{\infty}(_{\lambda}\X, \V) )$ when $\lambda'\uparrow \lambda$. This proves the left-continuity at~$\lambda>\lambda_c$.
\end{proof}

\subsection{Robust percolation in the limit model}

In this subsection, we work in the supercritical phase $\gamma >
\frac{\delta}{1+ \delta}$. We want to show that the infinite graph $G^\infty$, as well as its thinned versions $\tensor[^p]{G}{^\infty}$, for any $p>0$, do contain an infinite component.
\medskip

We sketch a simple strategy to get an infinite path in  $G^\infty$. Start from a sufficiently old vertex. Use a vertex born after time $\frac12$ as a stepping stone
to connect the old vertex by two edges to a much older vertex. Keep going forever, moving to older and older vertices. To ensure that this procedure generates an infinite path with 
positive probability we need to show that an old vertex is \wep at graph distance two in $G^\infty$ from a much older vertex, so that the failure probabilities sum to a value
strictly less than one. To get the necessary estimates we also need to avoid using old vertices that have an exceptionally small degree. The expected degree of a vertex with birth 
time~$s$ is of order~$s^{-\gamma}$, and the lemma below allows to define a notion of good vertices in such a way that (i)~every good vertex satisfies a lower bound on the degree at time~$\frac12$,
and (ii)~most old vertices are good. To this end fix a function $g$ as in Lemma~\ref{goodness_prob_bound} of the appendix, we will never need to know more about it than the fact it grows slower than 
polynomially.

\begin{defin}
For a vertex $\x = (x,s)$ in $G^\infty$ we denote by $\indeg ( \x ,t)$ its indegree at time $t\geq s$, or in other words the number of vertices born during the time interval $[s,t]$ that are connected by an edge to~$\x$.  The vertex $\x = (x,s)$ is a  \textbf{good vertex} if $s<1/2$ and 
$$\indeg ( \x ,1/2) \ge s^{- \gamma}/ g(s^{-1}).$$
\end{defin}


Lemma~\ref{goodness_prob_bound} ensures that under $\P_x$ the vertex~$\x$ with birth time~$s\le \frac12$  is good $\whp(1/s)$. 

\begin{defin}
 A vertex $\x = (x,s)$ born at time $s<1/2$ is \textbf{locally good} if its indegree in the graph $G^\infty(\X \cap\left([x-s^{-1/d}, x+s^{-1/d}]^d\times (0,1/2]\right), \V)$ is at least equal to $s^{- \gamma}/ g(s^{-1})$.
\end{defin}

\pagebreak[3]

The advantage of this more restrictive definition is that we can ensure that a vertex is locally good by only watching for the set of vertices nearby, up to distance $s^{-1/d}$. Moreover, by Lemma~\ref{goodness_prob_bound}, a vertex $(x,s)$ is still locally good $\whp(1/s)$.
The next lemma quantifies how young vertices allow to connect old vertices that have reached a high degree. 

\begin{lemma}[Two-connection lemma]\label{jumptohub}
Let $\x_{1}$ and $\x_{2}$ be two vertices of $G^\infty$ born before time~$1/2$.   
Write $\ell:=\tmtextsf{d} ( \x_{1} , \x_{2} )^{d}$ and $z_{i}:=\indeg (\x_{i} ,1/2)$
assuming $\ell= \Omega (z_{i} )$, for $i\in\{1,2\}$. Suppose there exists $\eps >0$ such that
$$z_{1} z_{2}^{\delta} \ell^{- \delta} = \Omega \big(z_{1}^{\eps} \big).$$
Conditional on the restriction of $G^\infty$ to vertices born before time $1/2$, {\wep}$(z_{1} )$, 
the vertices $\x_{1}$ and $\x_{2}$ are connected through a vertex born after time $1/2$, in the 
graph $G^\infty$. The analogous result holds also for the thinned graphs $\tensor[^{p}]{G}{^\infty}$, for any given $p>0$.
\end{lemma}

\begin{proof}
  From the construction rule and the hypothesis $\indeg ( \x_{1} ,1/2)
  =z_{1}$, it is clear that every vertex $\x = (x,s)$ of $\X$ satisfying the
  conditions $s>1/2$, $\tmtextsf{d} ( \x , \x_{1} )^{d} \le f (z_{1} )$ and
  $\V_{\x , \x_{1}} \le \varphi (1)$ is connected to $\x_{1}$. The number of
  such vertices is a Poisson variable with parameter of order~$z_{1}$,
  therefore it is {\wep}$(z_{1} )$ of order $z_{1}$.
The probability that $\x$ connects to $\x_{2}$ is
  \[ \varphi \left( \frac{s \tmtextsf{d} ( \x , \x_{2} )^{d}}{f( \indeg (
     \x_{2} ,s-))} \right) \ge \varphi \Big( \frac{(f(z_{1} )^{1/d} +\ell^{1/d}
     )^{d}}{f(z_{2} )} \Big) . \]
  On the right hand side, the numerator of the argument of $\varphi$ is $O (\ell)$, and the
  Potter bounds for regularly varying functions, see~\cite[Theorem~1.5.6]{bingham}, 
  ensure the right hand side is $\Omega ((\ell z_{2}^{-1} )^{- \delta - \eps'}
  )$, for any $\eps' >0$.
  Hence the event that one of the vertices~$\x$ is connected to~$\x_{2}$ by an edge is stochastically bounded from below
by a binomial random variable with parameters of order  $\Omega (z_{1} )$ (number of trials) and 
$\Omega ((\ell z_{2}^{-1} )^{- \delta -  \eps'} )$ (success probability). If $\eps'$ is chosen
  small enough, the expectation of this random variable is \smash{$\Omega(z_{1}^{_{\eps /2}} )$}, and therefore 
it is {\wep}$(z_{1} )$ of order \smash{$\Omega (z_{1}^{_{\varepsilon /2}} )$}. In particular, it is positive, and 
it stays {\wep}$( z_{1} )$ positive after percolation with any retention parameter~$p>0$.
\end{proof}


\begin{corollary}\label{cor_jumptohub}
  Suppose $\gamma > \frac{\delta}{1+ \delta}$. Choose 
$$\mbox{ first }  \alpha\in \big( 1 , \sfrac{\gamma}{\delta ( 1- \gamma )}\big) , \mbox{ then } \beta\in\big( \alpha , \sfrac{\gamma}{\delta}
( 1+ \alpha \delta )\big).$$
If the vertices $\x = ( x,s )$ and $\y = ( y,t )$ are good
  vertices with $t<s^{\alpha}$ and $\tmtextsf{d}( x,y )^{d} <s^{- \beta}$,  then they are connected through a vertex born
  after time $1/2$, {\wep}$( 1/s )$.
\end{corollary}

With this corollary in hand, we can state and prove the following key proposition.
\pagebreak[3]

\begin{proposition}
  [Chain of ancestors]\label{ancestors}

 If $\x_{0} = ( x_{0} ,s_{0} )$ is a locally good vertex, then
 \begin{enumerate}
 \item
   {\wep}$( 1/s_{0} )$ there exists a locally good vertex $\x_{1} = ( x_{1} ,s_{1} ) \in \X$ with $s_{1}
  <s_{0}^{\alpha}$ and $\tmtextsf{d}( x_{0} ,x_{1} )^{d} <s_{0}^{- \beta}$. We say $\x_{1}$ is an \emph{ancestor} of $\x_{0}$.\\[-3mm]

  \item  {\wep}$( 1/s_{0})$ there exists an infinite chain of ancestors $(\x_k)_{k\ge1}$, namely locally good vertices satisfying $s_{k+1}
  <s_{k}^{\alpha}$ and $\tmtextsf{d} ( x_{k} ,x_{k+1} )^{d} <s_{k}^{- \beta}$ for every $k\ge0$.\\[-3mm]

\item $\wep(1/s_0)$, two consecutive ancestors of the infinite chain of ancestors are always connected through a vertex born after time $1/2$, and therefore within graph distance two in $G^\infty$.
  \end{enumerate}

%
\end{proposition}

\pagebreak[3]

\begin{proof}
 The only difficult part is to explain how to find the ancestors. We have ensured that $\x_0$ is locally good, by looking only at the vertices in the set 
 $$\X\cap\Big([x_0-s_0^{1/d}, x_0+s_0^{1/d}]^d \times (0,1/2]\Big).$$ 
We always search for ancestors by moving to the right on the first coordinate. Let~$\eps=\beta-\alpha$. 
Take $\lfloor s_0^{_{-\eps/d}}/6\rfloor -1$ disjoint intervals of length $6s_0^{_{-\alpha/d}}$ inside \smash{$[s_0^{_{-\alpha/d}}, s_0^{_{-\beta/d}}]$}.  
Write $a_k$ for the centre of the $k$-th interval and $A_k:=(a_k, 0,\dots, 0)+[-3s_0^{_{-\alpha/d}},3 s_0^{_{-\alpha/d}}]^d$. The blocks $x+A_k$ are disjoint and 
have not been observed so far. Therefore, independently of everything else, the probability that $\X\cap(x+A_k)$ contains a locally good vertex at distance less 
than $s_0^{_{-\alpha/d}}$ from the centre of the block with birth time in $(s_0^\alpha/2, s_0^\alpha)$ is bounded from zero, say by~$c>0$.
One of the $\lfloor s_0^{_{-\eps/d}}/6\rfloor -1$ independent trials with success probability $\geq c$ has to succeed, $\wep(1/s_0)$, which proves~(1). Similarly, given the $k$ first ancestors, we find the $(k+1)$-th ancestor $\wep(1/s_k)$. Therefore we easily see that we can find an infinite chain of ancestors, $\wep(1/s_0)$, which proves~(2) as well.
\end{proof}

It follows directly from Proposition~\ref{ancestors} that the infinite graph
percolates in the supercritical phase. The same proof also holds for the thinned infinite graphs $\tensor[^p]{G}{^\infty}$, so that they also percolate. 
This immediately implies that $\tensor[^p]{\theta}{}>0$, for any $p>0$. Of course, the only infinite component of $\tensor[^p]{G}{^\infty}$ is a subgraph 
of that of $G^\infty$. In this sense, the infinite component of the infinite graph exists, is unique and robust.

\subsection{Robustness of the giant component} \label{finite_robustness}

We now show the following result.

\begin{proposition}\label{finerob}
Let \smash{$\gamma> \frac{\delta}{1+\delta}$} and $p>0$. With high probability, the largest component of \smash{$\tensor[^p]{G}{^t}$ contains $(\tensor[^p]{\theta}{}+ o(1))t$} vertices, 
while the second largest contains only $o(t)$ vertices. Hence there is a unique giant component, which has asymptotic density $\tensor[^p]{\theta}{}>0$.
\end{proposition}
 
It is tempting to believe that the result follows from the robustness of the infinite component of the infinite graph
by a pure approximation argument. 
However, the equivalence of the existence of a unique infinite cluster in the infinite graph, 
and that of a unique giant component in the finite graphs, does not always hold for long-range percolation models. In the rest of this subsection 
we prove this equivalence with a significant  effort, and our argument  uses specific features of the graphs in the phase \smash{$\gamma>\frac{\delta}{1+\delta}$}. 
The proof is organized in three parts. First, a direct study of the finite picture shows that the old vertices form a core of very-well connected vertices. The method is similar to that of the `chain of ancestors' argument for the limit model. Then, we explain how the law of large numbers allows to pass from the infinite to the finite picture. Finally, we have 
to show a convergence result for a carefully chosen functional to complete the proof.

\subsubsection*{Core of well-connected vertices}

With high probability, the oldest vertex of $G^t$ has birth time within \smash{$(\frac1{t \log t}, \frac{\log t}{t})$}. Then Lemma~\ref{goodness_prob_bound} ensures it is also good with high probability\footnote{The reader may note that the function $\log$ used here and in Lemma~\ref{goodness_prob_bound} plays no special role other than being a function growing to infinity slower than polynomially. And indeed, here and in Lemma~\ref{goodness_prob_bound}, we could have replaced it by any other function growing to infinity slower than polynomially.}. We now work implicitly on this event. In particular, we say that a statement holds \wep while the precise formulation should be `\wep, the statement holds or the oldest vertex is not a good vertex with birth time within 
\smash{$(\frac1{t \log t}, \frac{\log t}{t})$}'. Such an abuse of notation is used in the following proposition. For its formulation define, for $k>0$, the \emph{$2k$-core} to be the set of good vertices~$\x$ with birth time \smash{$s<t^{-1/\alpha^k}$.}
\pagebreak[3]

\begin{proposition} \label{core}\ \\[-5mm]
\begin{itemize}
\item[(1)]  In $G^{t}$, {\wep}, every vertex of the $2k$-core is connected, through a vertex born after time~$1/2$, 
to an ancestor or to the oldest vertex (or is the oldest vertex). \\[-3mm]
\item[(2)] In $G^t$, \wep, every vertex that has reached degree at least $t^{\gamma/\alpha^k}$ at time $1/2$ is connected  through a vertex born after time
  $1/2$ to the $2k$-core.
\end{itemize}
\end{proposition}

\begin{proof}
For item~(1), proceed as in Proposition~\ref{ancestors} to ensure that, $\wep(1/s)$, a good vertex $( x,s )$ with birth time \smash{$s<t^{-1/\alpha^k}$} has an ancestor and is connected to it through a vertex born after time~$1/2$. The slight difference is that we work here in the finite graphs and on the torus. The proof has to be adapted when $s^{-\beta/d}>\frac12\, t^{-1/d}$, because then the blocks we consider cover the whole torus and overlap.
In that case, recall that the oldest vertex is a good vertex, within distance $\frac12\, t^{1/d}$ of $x$, and born before time \smash{$\frac{\log t}{t}$}.
Corollary~\ref{cor_jumptohub} then ensures it is \wep connected, through a vertex born after time $1/2$, to $(x,s)$.
Finally, item~(2) follows easily from a further application of Lemma~\ref{jumptohub}.
\end{proof}

Proposition~\ref{core} has an important consequence, namely, \wep, all the vertices of the $2k$-core, as well as all the vertices that have reached at least  degree $t^{\gamma/ \alpha^k}$ at time 1/2, belong to the same connected component of $G^{t}$. Moreover, any two such vertices are within distance $4k+4$.
The connected component of the core is a natural candidate for the
giant component of the graph $G^{t}$. 

\subsubsection*{Use of the laws of large numbers}
Recall that Formula~\eqref{finitecomponents} gives, \whp, the asymptotic density of vertices of $\tensor[^p]{G}{^t}$ belonging to components of size $<k$, for any given $k\in \N$. As $k\uparrow\infty$ it goes to $1-\tensor[^p]{\theta}{}$ and hence, \whp an asymptotic proportion $(1-\tensor[^p]{\theta}{})$ of vertices belongs to finite-size components. The remaining $\tensor[^p]{\theta}{} t$ vertices belong to large clusters, whose sizes grow with $t$. However, at this point, nothing guarantees that they form one giant component. They could belong to various components of logarithmic size, for example.
\medskip

To see why this is not the case, we search for an indicator function $\xi_t$ taking the value one on exactly one component of $\tensor[^p]{G}{^t}$, which converges in probability to $\xi$, the indicator function of 
the event that $\x$ belongs to the infinite component of $\tensor[^p]{G}{^\infty}$. Inspired by the description of the core, we define $\xi_t$ to be the indicator of the event that $\x$ is connected through a path to the oldest vertex of $\tensor[^p]{G}{^t}$. If we can prove the convergence in probability of $\xi_t((0,U), \tensor[^p]{G}{^t})$ to $\xi((0,U),\tensor[^p]{G}{^\infty})$, then the law of large numbers gives
\smash{$ \frac 1 {pt} \sum_{\x \in \tensor[^p]{\X}{}^t} \xi_t(\x, \tensor[^p]{G}{^t}) \rightarrow \tensor[^p]{\theta}{}.$} In other words, \whp the component of the oldest vertex contains $(\tensor[^p]{\theta}{}+o(1))t$ vertices, as claimed. 
By the previous paragraph the second largest component cannot contain an asymptotically  positive proportion of all the vertices. This completes the proof of Proposition~\ref{finerob} subject to the assumed
convergence, which we now prove. In the proof we assume  $p=1$ to lighten the notation; the general case
follows by the same line of arguments. 
\pagebreak[3]

\subsubsection*{Convergence of $\xi_t$ to $\xi$}
Recall the notation $\Pp(du \, d\omega)=\int_0^1 du\, \P_{(0,u)}(d\omega)$ and $U(u,\omega)=u$. We have to prove the following two statements:
$$
\Pp ( A_{t} \backslash A_{\infty} )   \underset {t\to \infty} \longrightarrow  0 \quad \text{ and }\quad
\Pp ( A_{\infty}\backslash A_{t} )   \underset {t\to \infty} \longrightarrow  0,
$$
where $A_t$ denotes the event that $(0,U)$ is connected to the oldest vertex of $G^{t}$, and $A_\infty$ denotes the event that $(0,U)$ belongs to the infinite component of $G^\infty$.
The first statement follows directly from the almost sure local convergence of $G^t$ to $G^\infty$. Indeed, on the complement of the event $A_{\infty}$, the
component of $( 0,U )$ in $G^\infty$ is finite. For large enough $t$, it coincides with its component in the finite picture. Increasing $t$ further, if necessary,  we get that  this component does not contain the oldest vertex. 
\medskip

\pagebreak[3]

The second statement is significantly harder to prove. We first show that on $A_{\infty}$, the vertex $(0,U)$ is still in the infinite component of the graph $_\lambda G^\infty$, introduced in Subsection~5.2, if $\lambda$ is slightly less than 1. The remaining vertices in $\X\backslash$ $\!\!_{\lambda}\X$ will be used at the end of the proof for a sprinkling argument.
\medskip

Recall that $_\lambda G^\infty$ has the same law as $G^\infty$ but with $f$ replaced by $\lambda f$, and thus $\gamma$ by~$\lambda \gamma$. Taking \smash{$1>\lambda>\frac{\delta}{\gamma(1+\delta)}$}
we thus get $_\lambda\theta>0$. In particular we infer that $\lambda_c<1$ and hence, by Proposition~\ref{con},  the mapping $\mu \mapsto$ $\!\!_\mu\theta$ is left-continuous at the point $\mu=1$. Denote by  $_\lambda A_\infty$ the event that $(0,U)$ is connected to infinity in the graph $_\lambda G^\infty$. Then left-continuity implies that $\Pp ( A_{\infty} \backslash$ $\!_\lambda A_{\infty} )\to 0$ as $\lambda\uparrow1$. Hence it suffices to show, for fixed $\lambda\in(\frac{\delta}{\gamma(1+\delta)},1)$, that
$$\Pp \big(\,  \!_\lambda A_{\infty} \backslash A_{t} \big)   \underset {t\to \infty} \longrightarrow  0.$$
The remainder of the proof being technical and geometric, we write it for the case~$d=1$, for the sake of clarity. No argument is specific to dimension one, however, and it is not hard to see that the proof works, \emph{mutatis mutandis,} in higher dimensions.
\medskip

We fix a small parameter $a\in(0,1)$ and a large parameter $m\in(1,\infty)$. Precise constraints on these parameters will be given later. We show that on $_\lambda A_\infty$, the 
vertex $(0,U)$ is likely to be connected in the finite graph $_\lambda G^{t^m}$ to a \emph{moderately old} vertex born before time $t^{-a}$, and then use the sprinkling argument to ensure that this vertex is likely to be connected to the core, and thus to the oldest vertex of $G^{t^m}$.
Let $_\lambda B_t$ be the event that neither in $_\lambda G^{t^m}$ nor in $_\lambda G^\infty$
a vertex located in $[-t,t]$ is incident to an edge of length larger than $t^m/2-t$. 
On the event $\{U> t^{-a}\}$, which holds \whp, define $_\lambda C_0$, the component of $(0,U)$ in
the subgraph of $_\lambda G^\infty$ whose vertex set is restricted to
$_\lambda\X^t \cap \left([-t,t]\times (t^{-a}, 1]\right)$. Let $_\lambda\hat{A}_{t} \subset \{U>t^{-a}\}$ be the event that $_\lambda C_0$ is connected 
in $_\lambda G^\infty$ by a direct edge to a vertex in $[-t,t]\times (0,t^{-a}).$ 
\medskip

The proof is carried out in three steps:
\begin{itemize}
\item {\bf First step:}  $\Pp \big(\!$ $\!_\lambda B_t^{\rm c} \big)   \underset {t\to \infty} \longrightarrow  0$.
\item {\bf Second step:}  $\Pp \big( (\!$ $\!_\lambda {A}_{\infty} \cap$ $\!_\lambda B_t )\backslash\!$  $_\lambda\hat{A}_{t} \big)   \underset {t\to \infty} \longrightarrow  0$.
\item {\bf Third step:}  $\Pp  \big( (\!\!$ $_\lambda\hat{A}_{t} \cap\!$ $_\lambda B_t)\backslash  A_{t^m}  \big)   \underset {t\to \infty} \longrightarrow  0$.
\end{itemize}

\begin{proof}[Proof of the first step.]
The set $_{\lambda}\X^t \cap \left([-t,t]\times (0, 1]\right)$ consists of a Poisson number of vertices, with mean~$2\lambda t$, all with birth time uniform in $(0,1)$. The probability that a vertex with birth time uniform in $(0,1)$ is incident to an edge of length larger than $K$ has been estimated in~\cite{JM14}, see Theorem 4 and its proof, and is bounded by a constant times $K^{-\eta_\lambda}$, where $\eta_\lambda$ is the smallest of the three constants $1$, \smash{$\frac1{\lambda \gamma} - 1$} and $\delta-1$. In the robust phase, $\eta_\lambda\in (0,1)$. Taking $K=t^m/2-t$, we easily see that $\whp$ no vertex in $_\lambda\X^t \cap \left([-t,t]\times (0, 1]\right)$ is incident to an edge of length larger than $t^m/2-t$, as soon as the constant $m$ is chosen larger than $1/\eta_\lambda$. This reasoning works in $_\lambda G^\infty$ and in $_\lambda G^{t^m}$ as well. This proves the first step, under the constraint~$m>1/\eta_\lambda$.
\end{proof}

\begin{proof}[Proof of the second step.]
As $\Pp(U\leq t^{-a})\to 0$, we suppose $U>t^{-a}$. We fix $b>0$, to be specified later, and split our event in two parts, depending
on the value of $\vert$ $\!\!_\lambda C_0 \vert$. 
\medskip

{\bf Part~(i):} $\Pp\big( \{\vert$ $\!\!_\lambda C_0\vert < t^b\}\cap$ $\!_\lambda{A}_{\infty} \cap$  $\!_\lambda B_t \backslash$  $\!_\lambda\hat{A}_{t} \big)   
\underset {t\to \infty} \longrightarrow  0.$\\[-1mm]

On $_\lambda{A}_{\infty}$, a vertex of $_\lambda C_0$ has to be connected by an edge to a vertex outside $[-t,t]\times(t^{-a},1]$. On $_\lambda B_t \backslash$
$\!_\lambda\hat{A}_{t}$, this vertex has to be located in $[-t^m/2,t^m/2]\backslash[-t,t]$. Hence one of the vertices of $_\lambda C_0$ is incident to an edge longer 
than $t^{1-b}$, called a \emph{long edge}. This long edge either links two vertices of $_\lambda C_0$, or one vertex of $_\lambda C_0$ to a vertex located in 
$[-t^m/2,t^m/2]\backslash[-t,t]$. It is easy to see that a given vertex is unlikely to be incident to a long edge. But we can also prove that among all the vertices of 
$[-t,t]\times(t^{-a},1],$ many are incident to long edges. Therefore in this proof we must use the fact that $_\lambda C_0$ has few vertices (no more than $t^b$) and 
check that these vertices 
are not incident to long edges. In order to do that, we explore and reveal the component $_\lambda C_0$ and 
control at each step the probability of finding a long edge.
\pagebreak[3]
\medskip

We first explain the exploration process and how the information about the graph is progressively revealed.
When exploring the neighbourhood of a vertex we use the term \emph{inedge} to denote an edge connecting the vertex to a younger neighbour. The position and birth times of all vertices are all revealed 
at once, that is, we work conditionally on $_\lambda \X$, and the remaining randomness is only encoded in the variables $(\V(\x,\x'))_{\x,\x'\in{} _\lambda \X}$. For $\x\ne \x'\in{} _\lambda \X$, 
the indegree evolution processes $(Z(\x, t))_{t\ge s}$ and $(Z(\x', t))_{t\ge s'}$ are conditionally independent. Start the exploration with the single vertex $(0,U)$. Reveal all 
its neighbours. If $(0,U)$ is incident to a long edge, then stop the exploration and declare you found a long edge. If it is not the case, $(0,U)$ is declared `explored', 
its neighbours are declared `to explore'. Now choose a vertex left to explore. Reveal all its neighbours, except that you do not reveal whether it is connected by an edge 
to an older vertex `to explore'. (That edge will be checked only when we explore the inedges of the older vertex). If you have revealed a long edge, stop. Otherwise, 
the new neighbours you have revealed are added to the set of vertices `to explore', and the vertex is declared explored. Continue until there are no vertices left to explore, 
or a long edge is found.
An important feature of this exploration process is that it will eventually reveal all of $_\lambda C_0$ in at most $|_\lambda C_0|<t^b$ steps, unless it has been stopped for finding 
a long edge. Moreover, at each step, the information gathered about the indegree evolution of the vertices is controlled in the following way. For each vertex that 
is neither explored nor to explore, we have revealed the absence of some inedges (those that could have linked it to an explored vertex, precisely). For each vertex 
left to explore, we have revealed the presence of \emph{at most one inedge}, and the absence of several other inedges. Now we have to bound the probability of finding 
a long edge, conditionally on this information.
\smallskip

We introduce the following notation,
\begin{eqnarray*}
\Y_1&=& _\lambda\X \cap \big( [-t,t]\times(t^{-a},1] \big), \\
\Y_2&=& _\lambda\X \cap \big(([-t^m/2,t^m/2]\backslash [-t,t]) \times(t^{-a},1] \big), \\
\Y_3&=& _\lambda\X \cap \big( ([-t^m/2,t^m/2]\backslash [-t,t])\times(0,t^{-a}) \big).
\end{eqnarray*}
Though the vertices of $_\lambda C_0$ that we explore are all in $\Y_1$, the vertices of $\Y_2$ and $\Y_3$ have to be considered as potential endvertices of long edges.
With high (and even extreme) probability, no vertex of $\Y_1\cup \Y_2$ has reached indegree more than $t^{a\gamma}$, and no vertex $(y,s_y)$ of $\Y_3$ has reached indegree 
more than $s_y^{-\gamma}$, so we work on this event.
Without extra conditioning, bounds on the connection probabilities are easy to establish.
Indeed, if  $\x,\y\in\Y_1\cup\Y_2$ and $\eps>0$, we can roughly bound the probability that they are connected by
\begin{equation}\label{longedge12}
P(\x\leftrightarrow \y) \le \varphi(t^{-a} \tmtextsf{d}(\x,\y) / f(t^{a \gamma}))
 \le c_\eps t^{a(1+\gamma)(\delta-\eps)} \tmtextsf{d}(\x,\y)^{-(\delta-\eps)},
\end{equation}
where $c_\eps$ is given by the Potter bounds~\cite[Theorem~1.5.6]{bingham}. Similarly, if $\x\in\Y_1$ and $\y=(y,s_y)\in \Y_3$, then
\begin{align}
\nonumber P(\x\leftrightarrow \y) &\le \varphi(t^{-a} \tmtextsf{d}(\x,\y) / f(s_y^{-\gamma}))
\le \left(c_\eps t^{a(\delta-\eps)}s_y^{-\gamma(\delta-\eps)} \tmtextsf{d}(\x,\y)^{-(\delta-\eps)}\right) \wedge 1\\
\label{longedge3} &\le c_\eps t^{a(\delta-\eps)}  \left(s_y^{-\gamma(\delta-\eps)} \tmtextsf{d}(\x,\y)^{-(\delta-\eps)} \wedge 1\right).
\end{align}
Now suppose $\x\in\Y_1$ is the vertex we are currently exploring in the exploration process, and $\y$ is a vertex at distance $\ge t^{1-b}$, whose connection to $\x$ we have to check.
If $\y \in \Y_3$, then $\y$ is the older vertex, and its indegree evolution process is only conditioned on the nonexistence of some edges. In that case it only decreases the probability that it is connected to $\x$, and we still can use the bound~\eqref{longedge3}. Similarly, if $\y\in \Y_1 \cup \Y_2$, and  $\y$ is the older vertex, or $\x$ is the older vertex but we have not revealed the  presence of any inedge of $\x$, then the conditional probability of $\x\leftrightarrow\y$ is still bounded by~\eqref{longedge12}.
\pagebreak[3]
\medskip

We give details only for the hardest case, when $\x$ has already an inedge revealed, and $\x$ is older than $\y$. We write $(\x,\y_0)$ for the inedge revealed, and $\y_1$, \ldots, $\y_n$ for the other vertices that have been revealed not to be linked to $\x$. We further condition on the values of $\V(\x,\x')$ for any $\x'$ different from $\y$, $\y_0$, $\y_1$, \ldots, $\y_n$, writing \smash{$\F'$} for the sigma-algebra generated by these random variables. Note that they do not determine whether $\x$ is linked to $\x'$ or not. However, if we know in addition that $\x$ is linked to $\y$, $\y_0$ and not to $\y_1$, \ldots, $\y_n$, then the inedges of $\x$ are all determined, as well as its indegree evolution process, which we write $z_+(\x, s)_{s\le s_x\le 1}$. If we know on the contrary that $\x$ is not linked to $\y$, nor to $\y_1$, \ldots, $\y_n$, and only to $\y_0$, then its deterministic indegree evolution process is written  $z_-(\x,s)_{s \le s_x \le 1}$.
The following computation is straightforward,
\begin{equation}
\label{ptotal} 
\begin{aligned}
P(\x\leftrightarrow \y, \x\leftrightarrow \y_0, \x \nleftrightarrow \y_1, \ldots, \x\nleftrightarrow \y_n \big\vert \F') &= p_+(\y) p_+(\y_0) \prod_{i=1}^n (1- p_+(\y_i)),\\
P(\x\nleftrightarrow \y, \x\leftrightarrow \y_0, \x \nleftrightarrow \y_1, \ldots, \x\nleftrightarrow \y_n \big\vert \F') &=(1-p_-(\y)) p_-(\y_0) \prod_{i=1}^n (1- p_-(\y_i)),
\end{aligned}
\end{equation}
where we have written $p_+(\y)$ for $\varphi(s_y \tmtextsf{d}(\x,\y) /f(z_+(s_y-)))$, namely the probability that $\x$ is linked to $\y$ knowing that its indegree evolution process has followed $z_+$ until then. Similarly, $p_-(\y) = \varphi(s_y \tmtextsf{d}(\x,\y) /f(z_-(s_y-)))$. There is an easy comparison between $z_+$ and $z_-$, namely
$$\begin{array}{lclr}
z_+(\x,s)&=&z_-(\x,s) & \text{if }s<s_y,\\
z_+(\x,s_y)&=&z_-(\x,s_y) +1 & \\
z_+(\x,s)&\ge&z_-(\x,s)+1 &\qquad \quad \text{if }s>s_y.
\end{array}$$
Hence, $p_+(\y)=p_-(\y)$, and for $i\in\{1,\ldots, n\}$, we have $(1-p_-(\y_i)) \ge (1- p_+(\y_i))$. Moreover,
$$ \frac {p_-(\y_0)} {p_+(\y_0)} = \frac {\varphi(s_{y_0} \tmtextsf{d}(\x,\y_0)/f(z_-(s_{y_0}-)) )}   {\varphi(s_{y_0} \tmtextsf{d}(\x,\y_0)/f(z_+(s_{y_0}-)) )}
\ge  \frac {\varphi(s_{y_0} \tmtextsf{d}(\x,\y_0)/f(0) )}   {\varphi(s_{y_0} \tmtextsf{d}(\x,\y_0)/f(t^{a \gamma}) )}
\ge  c'_\eps t^{-a\gamma(\delta+\eps)},$$
where the last inequality is ensured by the Potter bounds, with $c'_\eps$ a strictly positive constant depending~on $\eps$.
Now,
\begin{align*}
P(\x\leftrightarrow \y_0, \x \nleftrightarrow \y_1, & \ldots, \nleftrightarrow \y_n \big\vert \F')\\
& \ge \left(p_+(\y) + (1-p_+(\y)) c'_\eps t^{-a\gamma(\delta+\eps)}\right) p_+(\y_0) \prod_{i=1}^n (1- p_+(\y_i)) \\
&\ge c'_\eps t^{-a\gamma(\delta+\eps)} p_+(\y_0) \prod_{i=1}^n (1- p_+(\y_i)).
\end{align*}
Combining with~\eqref{ptotal}, we get
$P(\x \leftrightarrow \y \vert \x\leftrightarrow \y_0, \x \nleftrightarrow \y_1, \ldots, \nleftrightarrow \y_n, \F')\le
c_\eps'^{-1} t^{a \gamma (\delta+\eps)} p_+(\y).$
But $p_+(\y)$ is always bounded by~\eqref{longedge12}. Integrating with respect to the law of $\V_{\x,\x'}$ gives
\begin{equation}
\label{conditionallongedge12}
P\big(\x \leftrightarrow \y \big\vert \x\leftrightarrow \y_0, \x \nleftrightarrow \y_1, \ldots, \nleftrightarrow \y_n \big)
\le c_\eps c_\eps'^{-1} t^{a(2 \gamma\delta+\delta-\eps)} \tmtextsf{d}(\x, \y)^{-(\delta-\eps)}.
\end{equation}
Informally, the price to pay to have a bound for the conditional probability is at most the multiplicative factor
$c_\eps'^{-1} t^{a \gamma (\delta+\eps)}$.
%
Adding the inequalities on every $\y\in\Y_1\cup \Y_2\cup\Y_3$ such that $\tmtextsf{d}(\x,\y)\ge t^{1-b}$, we can bound the probability that $\x$ is incident to a long edge, conditionally on the beginning of the exploration process, by $E_1+E_2+E_3$, where
\begin{eqnarray*}
E_1&=& c_\eps c_\eps'^{-1} t^{a(2 \gamma\delta+\delta-\eps)} t^{-(1-b)(\delta-\eps)} \vert \Y_1 \vert,\\
E_2&=& c_\eps c_\eps'^{-1} t^{a(2 \gamma\delta+\delta-\eps)} \sum_{\y \in \Y_2} (\vert y \vert - t+t^{1-b})^{-(\delta-\eps)}, \\
E_3&=& c_\eps t^{a(\delta-\eps)} \sum_{\y \in \Y_3} \left(s_y^{-\gamma(\delta-\eps)} (\vert y \vert - t+t^{1-b})^{-(\delta-\eps)} \wedge 1\right).
\end{eqnarray*}
This bound is independent of $\x$. Hence the probability that the exploration process reveals a long edge in less than $t^b$ steps is bounded by $t^b(E_1+E_2+E_3)$. In other words, we have proved
that 
$$\Pp\big( \{\vert \, _\lambda C_0\vert < t^b\}\cap \, _\lambda{A}_{\infty} \cap
\,_\lambda B_t \backslash \, _\lambda\hat{A}_{t}\ \big\vert \, _\lambda\X \big) \le t^b(E_1+E_2+E_3).$$
In order to conclude~(i), we have to prove that the bound is likely to be small, that is, goes to zero in probability. As $\vert \Y_1\vert$ is \whp of order $t$, the first term $t^b E_1$ is \whp of order $t^{1-\delta+\eps+ a(2 \gamma\delta+\delta-\eps)+b(1+\delta-\eps)}$. If $a$, $b$ and $\eps$ are chosen small enough, this bound goes to zero. For the second and third one, we show their expectation goes to zero. We have
\begin{eqnarray*}
\E\Big[\sum_{\y \in \Y_2} (\vert y \vert - t+t^{1-b})^{-(\delta-\eps)}\Big] 
&=& 2 \lambda (1-t^{-a}) \int_{t}^{t^m/2} (y - t+t^{1-b})^{-(\delta-\eps)} dy,
\end{eqnarray*}
which is of order $t^{(1-b)(1-\delta+\eps)}$. Hence
$\E[t^b E_2]=O(t^{1-\delta+\eps+ a(2 \gamma\delta+\delta-\eps)+b(\delta-\eps)}),$
which also goes to zero if $a$, $b$ and $\eps$ are small enough. Finally,
\begin{align*}
\E\Big[ \sum_{\y \in \Y_3}  \big(s_y^{-\gamma(\delta-\eps)}  & (\vert y \vert - t+t^{1-b})^{-(\delta-\eps)} \wedge 1\big) \Big] \\
& = 2 \lambda \int_t^{t^m/2} dy \int_{0}^{t^{-a}} ds \, \big(s^{-\gamma(\delta-\eps)} (y - t+t^{1-b})^{-(\delta-\eps)} \wedge 1\big)\\
& =  2 \lambda \int_{t^{1-b}}^{t^m/2-t+t^{1-b}} dy \, \int_{0}^{t^{-a}} ds  \big(s^{-\gamma(\delta-\eps)} y^{-(\delta-\eps)} \wedge 1\big)\\
& \le  2 \lambda \int_{t^{1-b}}^{\infty} dy \, \Big(y^{-1/\gamma}+ \int_{y^{-1/\gamma}}^{1} s^{-\gamma(\delta-\eps)} y^{-(\delta-\eps)} ds \Big).
\end{align*}
Finishing the calculation, this bound is $O(t^{-\eta (1-b)})$ with $\eta=\min(1/\gamma-1,\delta-\eps-1)>0$
, and thus
$\E[t^b E_3]= O(t^{-\eta+ a(\delta-\eps) +b(1+\eta)})$ goes to zero if $a$ and $b$ are chosen small enough.
\bigskip

{\bf Part~(ii):} $ \Pp \big( ($ $\!\!_\lambda{A}_{\infty} \cap \{\vert$ $\!\!_\lambda C_0\vert\ge t^b\})\backslash$  $\!_\lambda\hat{A}_{t} \big)   
\underset {t\to \infty} \longrightarrow  0$.\\

On the event $_\lambda A_{\infty} \cap\{\vert _\lambda C_0\vert \ge t^b\}$, we work conditionally on $_\lambda C_0$ and try to connect each vertex of $_\lambda C_0$ to some vertex in $[-t,t]\times (0, t^{-a}]$. Fix some $\eps>0$. For any given vertex in $\X^t$, 
there exists \wep a vertex of  $_\lambda \X^t$ with birth time in $(0, t^{-a}]$ and within distance $t^{a+\eps}$, because the number of such vertices follows a Poisson law of parameter 
$2\lambda t^\eps$. 
For each vertex of $_\lambda C_0$, we may then choose such a vertex, and try to connect it, with success probability bounded below, independently of everything else, by $\varphi(t^{a+\eps}/f(0)).$ By the Potter bound, for each $\eps'>0$, this is bounded below by  $c_{\eps'} t^{-(\delta+\eps')(a+\eps)}$, where $c_{\eps'}$ depends only on $\eps'$.
The number of edges between $_\lambda C_0$ and $[-t,t]\times (0, t^{-a}]$ is therefore bounded from below by a binomial variable of parameters $\lfloor t^b\rfloor$ and $c_{\eps'} t^{-(\delta+\eps')(a+\eps)}$, hence it is positive \whp 
as soon as 
$b>(\delta+\eps')(a+\eps)$.
Reducing $a$ if necessary (as well as $\eps$ and $\eps'$), we can ensure that this inequality is satisfied, which concludes the second step.
\end{proof}
\pagebreak[3]
 
\begin{proof}[Proof of the third step.]
On $_\lambda\hat{A}_{t}$, the vertex $(0,U)$ is connected in $_\lambda G^\infty$ and within $[-t,t]$ to a vertex with birth time in~$(0,t^{-a})$. 
We choose arbitrarily such a vertex $\x=(x, s)$. On $_\lambda B_t$, all these connections remain in the finite graph $_\lambda G^{t^m}$, that is, $(0,U)$ is connected in this finite graph to~$\x$. It should be enough to say $\x$ is likely to be in the $2k$-core for a well-chosen~$k$, and thus connected to the oldest vertex of $\T_{t^m}$ by Proposition~\ref{core}. However, due to the complex way we used to find $\x$, it is not that easy to ensure it is a good vertex, or to say anything about its degree. This is where the sprinkling argument is used, and the reason why we have worked with $\lambda<1$ in the entire proof.\smallskip

We condition on $_\lambda G^{t^m}$ and on the choice of $\x=(x,s)$ with $s<t^{-a}$. The law of the graph
$$_{\lambda}\widetilde G^{t^m}:=G^{t^m}\big(\big(\X^{t^m} \backslash{_\lambda\X^{t^m}}\big) \cup\{\x\}, \V^{t^m}\big)$$
under this conditional law is also the (unconditioned) law of
$G^{t^m}($ $\!_{1-\lambda}\X^{t^m} \cup\{\x\}, \V^{t^m})$,
because the set $\X^{t^m} \backslash$ $\!_{\lambda}\X^{t^m}$ is a Poisson point process 
of intensity $1-\lambda$ independent of $_{\lambda}\X^{t^m}$.
As a consequence, we know that the vertex $\x$ has \whp reached degree at least $t^{a(1-\lambda) \gamma}/g(t^a)$ 
at time $1/2$ in $_{\lambda}\widetilde G^{t^m}$.
As $_\lambda G^{t^m}$ and $_\lambda\widetilde G^{t^m}$ are both
subgraphs of $G^{t^m}$, taking \smash{$k>\frac{\log(m/(1-\lambda)a)}{\log(\alpha)}$}
in Proposition~\ref{core} allows to conclude that $\x$ is \whp connected to the $2k$-core and 
in particular to the oldest vertex in $G^{t^m}$. Hence the same holds for $(0,U)$, and 
$A_{t^m}$ is satisfied.
\end{proof}
%

\section{Proof of non--robustness}

\subsection{Non--robustness for $\gamma<1/2$}

We have seen in Section~\ref{nonrob} that it suffices to show that $\tensor[^p]{G}{^\infty}$ contains no infinite component 
if $p$ is chosen small enough. We will introduce a notion of quick path, such that if $\tensor[^p]{G}{^\infty}$ contains an 
infinite component, then there exists an infinite quick path. Quick paths will be constructed in such a way that we can 
estimate their probability using a disjoint occurrence argument.

\subsubsection*{First moment method based on quick paths}
All the graphs we consider are locally finite, therefore an infinite component has to be of infinite diameter.
Actually, a vertex $\x_0$ of an infinite component is always the starting vertex of at least one infinite geodesic, 
that is an infinite path $(\x_n)_{n\ge 0}$, $\x_n=(x_n, s_n)$, in the graph with the property that
the graph distance between two vertices $\x_n$ and $\x_{n+k}$ is always $k$, for all $n,k \ge0$.
This can be proved with a simple diagonal argument that we leave to the reader.
Note that a geodesic is in particular a vertex and edge self-avoiding path.
Starting from any infinite geodesic $(\x_n)_{n\ge0}$ in $\tensor[^p]{G}{^\infty}$, 
we now construct deterministically, in two steps, an infinite self-avoiding path $(\z_n)_{n\ge0}$, $\z_n=(z_n,u_n)$, in $G^\infty$  called the quick path associated with $(\x_n)_{n\ge 0}$.
\medskip

First, we construct a subsequence $(\y_n)_{n\ge0}$ as $\y_n=\x_{\phi(n)}$. Start with $\phi(0)=0$, and thus $\y_0=\x_0$. Given 
$\phi(n)$, define
\begin{align*}
\skrin_n:=\big\{k> \phi(n) \colon  \exists\, \y=(y,t) \in G^\infty \mbox{ such that }  t>s_{\phi(n)}, s_k
\mbox{ and }  \x_{\phi(n)}{\leftrightarrow}\y{\leftrightarrow}\x_k \text{ in }G^\infty \big\}.
\end{align*}
The vertex $\y$ in the definition of $\skrin_n$ is called a \emph{common child} of the vertices $\x_{\phi(n)}$ and~$\x_k$, note that
it is chosen in $G^\infty$ and not in $\tensor[^p]{G}{^\infty}$.
For $k\in \skrin_n$, the graph distance between the vertices $\x_{\phi(n)}$ and $\x_k$ is thus (at most) 2 in $G^\infty$, while it is $\vert k-\phi(n)\vert$ in $\tensor[^p]{G}{^\infty}$. If $\skrin_n$ is non-empty, it has to be finite and we set $\phi(n+1)=\max \skrin_n$.
Otherwise, we set $\phi(n+1)=\phi(n)+1$.
\medskip

\pagebreak[3]

By its definition $(\y_n)_{n\ge0}$ satisfies the following properties:
\begin{itemize}
\item for all $n\ge0$ we have $\y_n\in \tensor[^p]{\X}{}.$
\item for all $n\ge0$ and $j\ge2$, the vertices $\y_n$ and $\y_{n+j}$ are not connected by an edge and have no common child in $G^\infty$.
\item for all $n\ge0$, the vertices $\y_n$ and $\y_{n+1}$ are either connected by an edge, or have a common child in $G^\infty$.
\end{itemize}
Finally, we create a third sequence by inserting, between every pair of vertices $\y_n$ and $\y_{n+1}$ that are not connected by an edge, the oldest common child in $G^\infty$. We obtain an infinite sequence \smash{$(\z_n)_{n\ge0}$}, which is an infinite self-avoiding path of $G^\infty$, and which we call the
quick path associated with $(\x_n)_{n\ge 0}$.
\medskip

We call  a vertex $\z_n$ in the quick path a \emph{regular vertex} if it is older than at least one of its neighbours $\z_{n-1}$ and $\z_{n+1}$, 
and we call $\z_n$ a \emph{local maximum} if it does not satisfying this property. Similarly define the \emph{local minima}.
Hence a vertex $\z_n$, with $n>0$, belongs to the sequence $(\y_k)_{k\ge0}$ if and only if it is regular.
With this terminology the path $(\z_n)_{n\ge0}$ has the following properties:

\begin{enumerate}
\item[(i)] Every regular vertex is  in $\tensor[^p]{\X}{}$. The starting vertex $\z_0$ is also in $\tensor[^p]{\X}{}$.
\item[(ii)] A regular vertex $\z_n$ cannot be connected by an edge to any younger vertex of the path, except possibly $\z_{n-1}$ and $\z_{n+1}$.
\item[(iii)] Two regular vertices $\z_n$ and $\z_{n+j}$, with $n\ge0$ and $j\ge 2$, can have common children only if $j=2$ and $\z_{n+1}$ is a local maximum.
In that case, $\z_{n+1}$ is their oldest common~child.
\end{enumerate}
Properties (ii) and (iii) depend only on the graph $G^\infty$ and not on the percolation procedure, and define the notion of 
\emph{quick paths}.
We also define the notion of quick paths for finite paths, by restricting the quantifiers accordingly.
Based on the observation that if $(0,U)$ is in the infinite component of~$\tensor[^p]{G}{^\infty}$ then it must be the starting point of arbitrarily long quick paths satisfying property~(i), the first moment calculation in the next subsection shows that the expected number of such paths of length $n$ goes to 0,
as $n\uparrow\infty$, if $p$ is small enough. Note that, given a quick path $\z_0{\leftrightarrow}\cdots{\leftrightarrow}\z_n$ in $G^\infty$, it satisfies condition~(i) with probability at most $p^{n/2}$, because at least $\frac{n}2-1$ of the vertices on a quick path of length~$n$ are regular. Therefore, if we show that the expected number of quick paths of length $n$ grows at most exponentially, namely it is $O(C^n)$ for some finite constant $C$, we can infer from Borel-Cantelli that, for any $p<1/C^2$, almost surely the component of~$(0,U)$ in $\tensor[^p]{G}{^\infty}$ is finite.

\subsubsection*{The expected number of quick paths has at most exponential growth.}

The expected number of quick paths of length $n$ is given by the multiple integral
\begin{align*}
\int_0^1 \E_{(0,u)}\big[ & \big\vert\{\text{quick paths of length }n\text{ starting at }(0,u)\} \big\vert\big] \, du\\
&= \int_{(0,1]\times (\R^d\times (0,1])^n} du_0 \, \, d\z_1 \ldots d\z_n \,
\P_{\z_0,\ldots,\z_n}\big( (\z_k)_{0\le k \le n} \text{ is a quick path in }G^\infty\big),
\end{align*}
where $\z_0=(0,u_0)$ and $\P_{\z_0,\ldots,\z_n}$ is the measure $\P$ conditioned on the event $\{\z_0,\ldots,\z_n \in \X\}$.
Under this measure, $\X$ is simply a Poisson point process of intensity one with the points $\z_0,\ldots,\z_n$ added.
The aim is now to bound the probability of $(\z_k)_{0\le k \le n}$ being a quick path, and see how we can integrate this bound.
Writing $\z_k=(z_k,u_k)$ for $0\le k\le n$, we will actually first integrate over $(z_1,\ldots,z_n)$ (space integration), then over $(u_0,\ldots,u_n)$ (time integration).
The main step will be to prove the following proposition.

\begin{proposition} \label{time_path_bound_prop}
There exists a finite constant $C$ such that for every $n\ge0$ and distinct numbers $u_0,\ldots,u_n$ in $(0,1]$, the following inequality holds,
\begin{equation} \label{time_path_bound}
\begin{aligned}
\int_{(\R^d)^n} dz_1 \ldots dz_n\,
\P_{(0,u_0), (z_1,u_1),\ldots,(z_n,u_n)} & \big( (z_k,u_k))_{0\le k \le n} \text{ is a quick path in }G^\infty\big) \\
&\le C^n \prod_{k=1}^n  \frac 1 {(u_{k-1}\wedge u_k)^\gamma (u_{k-1}\vee u_k)^{1-\gamma} }.
\end{aligned}
\end{equation}
\end{proposition}
The bound given in~\eqref{time_path_bound} is good in many respects. First, 
the proof provides a constant~$C$ that does not depend on the choice of the 
profile function~$\phi$, but only on the attachment rule~$f$.
Second, the term $1 /{(u_{k-1}\wedge u_k)^\gamma (u_{k-1}\vee u_k)^{1-\gamma}}$ is comparable to the probability that two vertices in 
a non-spatial equivalent of our model, with birth times $u_{k-1}$ and $u_k$, are connected by an edge.
But most importantly, the next lemma shows that after integrating over time we obtain the desired bound.

\begin{lemma} \label{time_integral_bound}
If $\gamma<1/2$, then there exists a finite constant $C'$ such that, for any $n>0$,
$$ \int_0^1du_0 \cdots\int_0^1  du_n \, \prod_{k=1}^n  \frac 1 {(u_{k-1}\wedge u_k)^\gamma (u_{k-1}\vee u_k)^{1-\gamma} } \le C'^n.$$
\end{lemma}

\begin{proof}
Pick $-\gamma>\alpha>\gamma-1$. Then carrying out the integration over $u_n$ gives
\begin{align*}\int_0^1du_0 \cdots\int_0^1  du_n \, & \, \prod_{k=1}^n  \frac {u_n^\alpha} {(u_{k-1}\wedge u_k)^\gamma (u_{k-1}\vee u_k)^{1-\gamma} }\\ 
& \le C'\int_0^1du_0 \cdots\int_0^1   du_{n-1} \, \prod_{k=1}^{n-1}  \frac {u_{n-1}^\alpha} {(u_{k-1}\wedge u_k)^\gamma (u_{k-1}\vee u_k)^{1-\gamma} },
\end{align*}
for $C'=(1+\alpha-\gamma)^{-1}-(\alpha+\gamma)^{-1}$. The result follows from this by induction.
\end{proof}

The proposition and the lemma, combined, prove the non-robustness of $(G^t)_{t>0}$ for any $\gamma<1/2$.
More precisely, if $C$ and $C'$ are the constants given in Proposition~\ref{time_path_bound_prop} and Lemma~\ref{time_integral_bound}, 
then for any $p<(C C')^{-2}$,  the graph $\tensor[^p]{G}{^\infty}$ contains no infinite component almost surely and, with high probability, 
the network $(\tensor[^p]{G}{^t})_{t>0}$ contains no giant component.

\subsection{Proof of Proposition~\ref{time_path_bound_prop}}

The proof of Proposition~\ref{time_path_bound_prop} is  based on two ingredients.
First, the definition of a quick path allows the use of a BK inequality, that splits the paths into small parts that interact with
negative correlation. Each small part comprises no more than four edges, and the probability of such a path can be bounded by more or less straightforward integration.


\subsubsection*{Splitting procedure}

We now explain how to split the sequence $\z_0, \ldots, \z_n$ into small parts.
The splitting procedure depends only on their birth times $u_0,\ldots,u_n$, which we assume to be pairwise distinct, 
but not on the spatial positions $z_0,\ldots,z_n$. The rule is simple, namely
$$
\begin{minipage}{.75 \textwidth}
\emph{For $i=0,\ldots,n,$ introduce a splitting at index $i$ if either $u_i$ is larger than both $u_{i-1}$ and $u_{i-2}$, or $u_i$ is larger than both $u_{i+1}$ and $u_{i+2}$.}
\end{minipage}
$$
\vspace{-2mm}

The boundary convention is that no condition is requested for indices outside $\{0,\ldots,n\}$, so that, for example,  there is always a splitting at indices 0 and $n$, and there is a splitting at index $1$ if $u_1>u_0$. We write $n_0=0< n_1<\cdots< n_k=n$ for the splitting indices in increasing order. These split the path into $k$ parts. The $j$th part consists of the sequence $(\z_{n_{j-1}}, \ldots, \z_{n_j})$, with $\z_{n_{j-1}}$ and $\z_{n_j}$ constituting the two \emph{boundary vertices} and the other vertices the \emph{inside} of part $j$. Note that the boundary vertices $\z_{n_1},\ldots, \z_{n_{k-1}}$ belong to two consecutive parts. A vertex $\z_i$ is a \emph{local maximum} of a part $(\z_{n_{j-1}}, \ldots, \z_{n_j})$ if  $n_{j-1}\leq i \leq n_j$ and we have both $u_i>u_{i-1}$ (if $i>n_{j-1}$) and $u_i>u_{i+1}$ (if $i<n_j$).   Observe that a boundary vertex of a part can be a local maximum of a part without being a local maximum. We say a vertex $\z_i$ \emph{contributes} to a part, if it belongs to, but is not a local maximum of, this part. 
\medskip

\pagebreak[3]

Using this terminology, we observe that
\begin{itemize}
\item Local maxima never contribute to any part (irrespective whether they are inside of a part or boundary vertices of two parts).
\item Local minima are always inside a part, and contribute to it.
\item The other vertices always contribute to exactly one part, whether they are inside it or at its boundary.
\end{itemize}
For $1\leq j \leq k$, let \smash{$A_j=\{\z_{n_{j-1}}{\leftrightarrow}\cdots{\leftrightarrow}\z_{n_j}\}$} 
be the event that \smash{$(\z_{n_{j-1}}, \ldots, \z_{n_j})$} is a path in~$G^\infty$.
Recall that a vertex $\x=(x,s)$ is a \emph{child} of $\z_i$ if $\z_i{\leftrightarrow}\x$ and $u_i<s$. We define~$\Sigma_j$ to be the (random) set of all children of vertices $\z_i$ contributing to part $j$, different from $\z_{i-1}$ and $\z_{i+1}$, 
which have birth times in the interval \smash{$(u_i, u_{i-1}\vee u_{i+1})$}.
Informally, the set $\Sigma_j$ contains all the information beyond  the variables $\V(\z_{i}, \z_{i+1})$, $n_{j-1}\leq i \leq n_j-1$, 
which is needed to check whether $A_j$ occurs or not.
\medskip\pagebreak[3]

The following lemma justifies the splitting rule.
\begin{lemma}\label{disjointness}
If $(\z_0,\ldots, \z_n)$ is a quick path, then the sets $\{\z_0,\ldots, \z_n\}$ and $\Sigma_1,\ldots, \Sigma_k$ are pairwise 
disjoint.
\end{lemma}
\begin{proof}
Property (ii) of the definition of quick paths implies that vertices of the path do not belong to any $\Sigma_j$. We now use
Property (iii) and the splitting procedure to see that $\Sigma_j$ and $\Sigma_{j'}$ do not intersect if $j{\neq}j'$. Indeed,
if they intersect, this would mean that a vertex~$\z_i$ contributing to part $j$ and a vertex~$\z_{i'}$ contributing to part $j'$ 
have a common child in $(u_i, u_{i-1}\vee u_{i+1})\cap (u_{i'}, u_{i'-1}\vee u_{i'+1})$. By Property~(iii) we must have $\vert i'-i\vert\le 2$.
If $\vert i'-i\vert=2$, say $i'=i+2$, then their oldest common child has to be $\z_{i+1}$.
As they contribute to different parts containing $\z_{i+1}$,  there must be a splitting at index $i+1$.
Hence either $u_{i-1}\vee u_{i+1}=u_{i+1}$, or $u_{i+1}\vee u_{i+3}=u_{i+1}$.
In each case, their common child with birth time in $(u_i,u_{i-1}\vee u_{i+1})\cap (u_{i+2},u_{i+1}\vee u_{i+3})$ is older than $\z_{i+1}$, and we get a contradiction. If $\vert i'-i\vert=1$, we can assume that $i'=i+1$ and there is a splitting at index $i$. 
The vertex~$\z_i$ cannot be a local maximum, as otherwise it contributes to no part.
Combining these two facts  we get that $(u_i, u_{i-1}\vee u_{i+1}) \cap (u_{i+1}, u_i\vee u_{i+2})$ is empty and hence a contradiction.
\end{proof}

\pagebreak[3]

\subsubsection*{BK--inequality}
We now use a version of the famous van den Berg-Kesten~(BK)~inequality to show that the probability of observing 
a quick path is bounded above by the product of the probabilities of the events $A_j$, for $j\in\{1,\ldots,k\}$. The version
of the BK--inequality  we use is valid for marked Poisson processes with unit intensity on a bounded domain, see~\cite{vandenberg}.
It states that the probability for increasing events $E_1,\ldots,E_k$ to occur disjointly is bounded above by the product of their individual probabilities, namely
$$ P(E_1\circ \ldots\circ E_k) \le \prod_{j=1}^k P(E_j).$$
In the present context, an event $E$ in the space of configurations of the marked Poisson processes
is called \emph{increasing} if, given any configuration in~$E$, the configuration  
with an arbitrary marked point added is also in~$E$. \emph{Disjoint occurrence} of events $E_1, \ldots, E_k$ is written 
as $E_1\circ \ldots\circ E_k$ and defined as  follows. 
A configuration consisting of the point set~$S$  with marks $(m_s)_{s\in S}$ is in $E_1\circ \ldots\circ E_k$, if we can find 
$S_1,\ldots,S_k$, disjoint subsets of $S$, so that, for each $j$, the set $S_j$ with marks $(m_s)_{s\in S_j}$ is in $E_j$.
We say the marked set $S_j$ ensures that $E_j$ is realized.
\medskip

Let us see how this inequality fits in our context. Recall that we work under $\P_{\z_0,\ldots,\z_n}$ and all the events we consider depend only on  the set of children of vertices $\z_0,\ldots,\z_n$. Therefore, the events are all deterministic functionals of the following ingredients:
\begin{enumerate}
\item the random variables $\V(\z_i,\z_{i'})$, for distinct indices $i, i' \in \{0,\ldots,n\}$;
\item the random set $\X':=\X\backslash\{\z_0,\ldots,\z_n\}$, which is a Poisson point process of unit intensity on $\R^d\times(0,1]$, together with the random marks $(\V(\x,\z_0),\ldots,\V(\x,\z_n))$ in $[0,1]^{n+1}$ attached to the vertices  $\x\in\X'$.
\end{enumerate}
\emph{First}, in order to study our problem in the framework of disjoint occurrence,  we have to remove the dependence on the random variables \smash{$\V(\z_i,\z_{i'})$}. We introduce \smash{$\P^{\ssup 0}_{\z_0,\ldots,\z_n}$} for the conditional probability given \smash{$\V(\z_i,\z_{i'})=0$}, for all $\vert i'-i\vert\ge2$. In other words, under \smash{$\P^{\ssup 0}_{\z_0,\ldots,\z_n}$} the indegree evolution process of a vertex $\z_i$ cannot grow because of vertices $\z_{i'}$ with $\vert i'-i\vert\ge2$, just as if the vertex $\z_i$ was not seeing them. We observe that
$$\P_{\z_0,\ldots,\z_n}\big((\z_n)\text{ is a quick path}\big) 
\le \P^{\ssup 0}_{\z_0,\ldots,\z_n}\big((\z_n)\text{ is a quick path}\big),$$
because if $(\z_n)$ is a quick path, then reducing the value of  $\V(\z_i,\z_{i'})$ with 
$\vert i'-i\vert\ge2$ neither affects properties~(ii) and~(iii), nor does it remove edges from the quick path.
We also work conditionally on the values $\V(\z_i,\z_{i+1})$, for $i\in\{0,\ldots,n-1\}$.%
\medskip

\emph{Second}, to apply the result of~\cite{vandenberg} we need to make the underlying domain bounded.
To this end  we work with the natural finite picture approximation of the graph and of our events.
For $t$ finite, but large enough so that \smash{$\T_t$} contains the vertices $\z_0,\ldots,\z_n$, we construct the graph $G^t$ and we denote by $A_j^t$ the event that the $j$th part is a path, in the graph $G^t$.
\medskip

Now the events $A_1^t, \ldots,A_k^t$ are increasing events of a marked Poisson point process with unit intensity on $\T_t\times(0,1]$. Applying the BK--inequality gives
$$\P^{\ssup 0}_{\z_0,\ldots,\z_n}\big(A_1^t\circ \ldots\circ A_k^t\, \big| \, 
\V_{\z_0,\z_1},\ldots, \V_{\z_{n-1},\z_n} \big)
\le \prod_{j=1}^k \P^{\ssup 0}_{\z_0,\ldots,\z_n} \big( A_j^t \, \big| \, \V_{\z_0,\z_1},\ldots,
\V_{\z_{n-1},\z_n} \big).$$
As $t\uparrow\infty$, we know that $G^t$ converges locally to $G^\infty$, and thus the indicator of 
\smash{$A_j^t$}  to that of ${A_j}$, almost surely.
Moreover, similarly to $\Sigma_j$, we define \smash{$\Sigma_j^t$} as the set of all children in $G^t$ of vertices $\z_i$ contributing to part $j$, different from $\z_{i-1}$ and $\z_{i+1}$, and with birth times in $(u_i, u_{i-1}\vee u_{i+1})$, then we also have that, almost surely, the sets $\Sigma_j^t$ coincide with the sets $\Sigma_j$, for $t$ large enough.
If $\z_0,\ldots,\z_n$ is a quick path of $G^\infty$, it is clear that for $t$ large enough, not only all the \smash{$A_j^t$} are satisfied, but also the sets \smash{$\Sigma_j^t$}, which ensure the events 
\smash{$A_j^t$} are satisfied, are disjoint, using Lemma~\ref{disjointness}.
Consequently, the events \smash{$A_j^t$} have to occur disjointly for $t$ large enough.
We get\\[-3mm]
$$\P^{\ssup 0}_{\z_0,\ldots,\z_n}\big(\z_0,\ldots,\z_n \text{ is quick}\, \big| \, 
\V_{\z_0,\z_1},\ldots,\V_{\z_{n-1},\z_n} \big)
\le \prod_{j=1}^k \P^{\ssup 0}_{\z_0,\ldots,\z_n} \big(A_j\, \big| \, 
\V_{\z_0,\z_1},\ldots, \V_{\z_{n-1},\z_n} \big).$$
The event $A_j$ depends actually only on $\V(\z_{n_{j-1}},\z_{n_{j-1}+1}), \ldots, \V(\z_{n_j-1},\z_{n_j})$.
Therefore, the product on the right hand side  is a product of independent 
random variables. Taking expectation, and using the tower property of conditional expectation, gives
$$\P^{\ssup 0}_{\z_0,\ldots,\z_n}\big(\z_0,\ldots,\z_n \text{ is a quick path} \big)
\le \prod_{j=1}^k \P^{\ssup 0}_{\z_0,\ldots,\z_n}(A_j).$$
Combining  with the observation that $\P^{\ssup 0}_{\z_0,\ldots,\z_n}(A_j)$ is also equal to $\P^{\ssup 0}_{\z_{n_{j-1}},\ldots,\z_{n_j}}(A_j)$, we obtain
\begin{equation}
\P_{\z_0,\ldots,\z_n}\big(\z_0,\ldots,\z_n \text{ is a quick path}\big)
\le \prod_{j=1}^k \P^{\ssup 0}_{\z_{n_{j-1}},\ldots,\z_{n_j}}(A_j).
\end{equation}

\begin{figure}[ht]
\begin{center}
\includegraphics[scale=0.5]{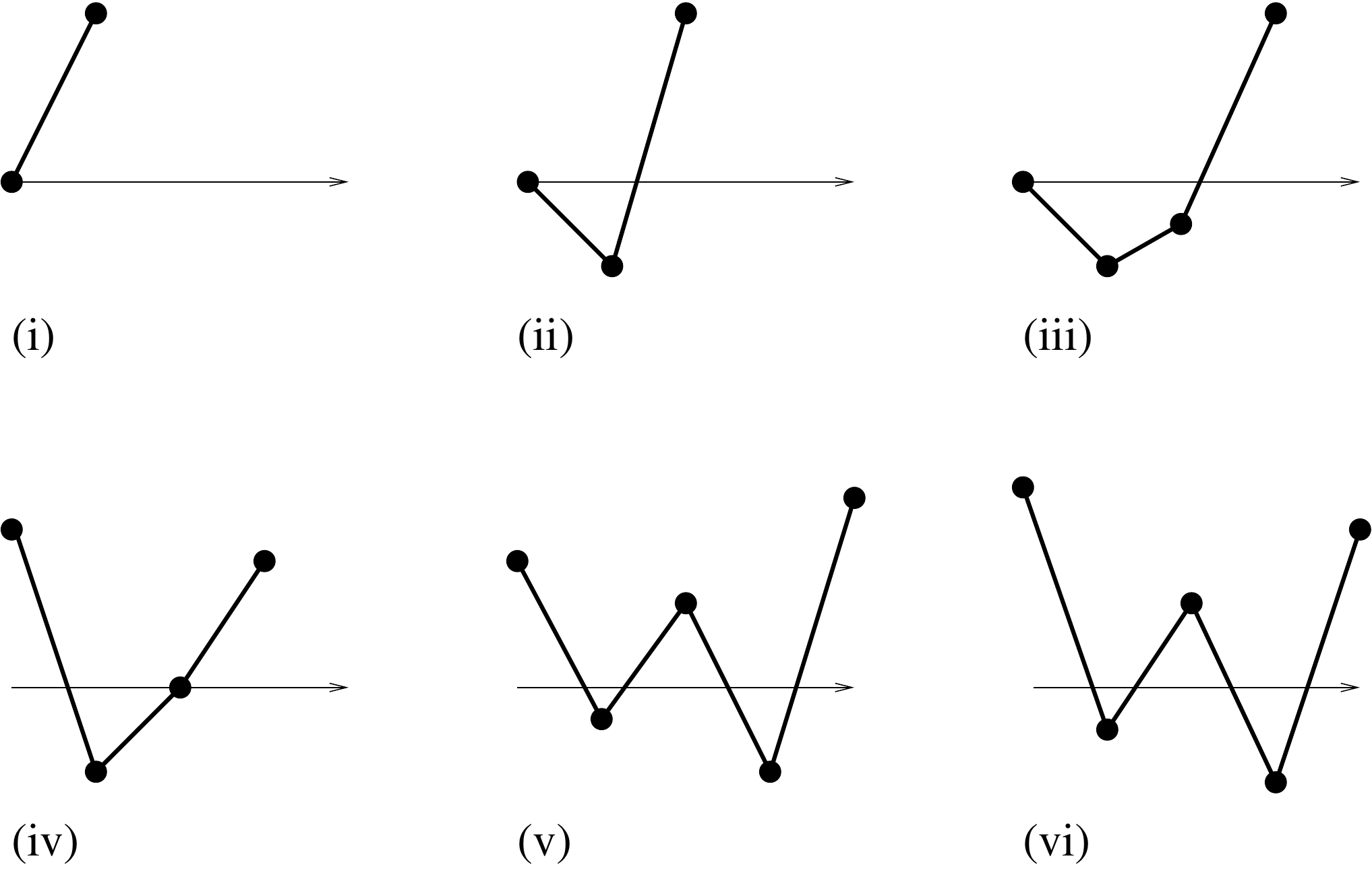}
\end{center}
\caption{Up to symmetry there are six types of small parts after the splitting. 
Illustrated, with the index of a point on the
abscissa and time on the ordinate, these are (i)~one single edge, 
(ii)~a {\sf V} shape with two edges,
(iii)~a {\sf V} shape with three edges and the end vertex of the short leg between the two vertices of the long leg,
(iv)~a {\sf V} shape with three edges and both vertices of the long leg below the end vertex of the short leg,
(v)~a {\sf W} shape with the higher end vertex on the side of the deeper valley,
(vi)~a {\sf W} shape with the lower end vertex on the side of the deeper valley.}
\end{figure}

\subsubsection*{Bound for the small parts}

We have bounded the probability of observing a (long) quick path by the product of the probabilities of observing a path, independently for each part.
In order to prove \eqref{time_path_bound}, it suffices to prove the corresponding inequality for each part $j$,
\begin{align*}
\int_{(\R^d)^{n_j-n_{j-1}}} dz_{n_{j-1}+1} \ldots dz_{n_j} \, &
\P^{\ssup 0}_{(z_{n_{j-1}},u_{n_{j-1}}),\ldots,(z_{n_j},u_{n_j})} 
\big( ((z_k,u_k))_{0\le k \le n} \text{ is a path in }G^\infty\big) \\
&\le C^{n_j-n_{j-1}} \prod_{k=1}^{n_j-n_{j-1}}  \frac 1 {(u_{k-1}\wedge u_k)^\gamma (u_{k-1}\vee u_k)^{1-\gamma} }.
\end{align*}
Instead of treating all six possible types of parts, listed in Figure~3, we only treat the most complex type,
numbered~(v). It should become clear to the reader that we are giving bounds in a way that 
would work similarly for all the other types. To lighten notation we also suppose the part is the first part,  that is, we suppose $u_3<u_1<u_2<u_0<u_4$, and show that
$$\int_{(\R^d)^4} dz_1 dz_2 dz_3 dz_4 \,
\P^{\ssup 0}_{\z_0,\z_1,\z_2,\z_3,\z_4} (A_1)
\le C^4 \frac 1 {u_0^\gamma u_1^{1-\gamma}} \frac 1 {u_2^\gamma u_1^{1-\gamma}} \frac 1 {u_2^\gamma u_3^{1-\gamma}} \frac 1 {u_4^\gamma u_3^{1-\gamma}}.$$
We introduce the canonical filtration $\F_t$, for $t\in(0,1]$, associated to the construction of the graph $G^\infty$ up to time $t$, i.e., $\F_t$ is the smallest $\sigma-$algebra for which the restriction of $G^\infty$ to vertices with birth times in $(0,t]$ is measurable. Similarly define $\F_{t-}$.
Observe that, writing $(Z_{\z}(t))_{t\in (0,1]}$ for the indegree evolution process of vertex $\z$, 
the process $(Z_{\z_1}, Z_{\z_3})$ is adapted to the filtration. In the following we use $c_1, c_2,\ldots$ to denote some positive constants depending only on the attachment rule~$f$.
\medskip

A change of variables from $(z_1,z_2,z_3,z_4)$ to $(y_1,y_1',y_3,y_3')=(z_2-z_1,z_0-z_1,z_2-z_3,z_4-z_3)$  and the tower property of conditional expectation yield
\begin{align*}
\int_{(\R^d)^4} dz_1 dz_2 & dz_3  dz_4 \, \P (A_1) \\
& = \int_{(\R^d)^4} dy_1 dy_1' dy_3 dy_3' \,  \E\big[ \E[\E[\E[\E[\E[\E[\one_{A_1} | \F_{u_4-}]|\F_{u_0}]|\F_{u_0-}]|\F_{u_2}]|\F_{u_2-}]|\F_{u_1}] \big],
\end{align*}
where we have here simply written $\E$ for expectation and conditional expectation under the probability measure \smash{$\P^{\ssup 0}_{\z_0,\z_1,\z_2,\z_3,\z_4}$}.
Rewriting the indicator of ${A_1}$ as product of indicators,  $\one{\{\z_0{\leftrightarrow}\z_1{\leftrightarrow} \z_2{\leftrightarrow} \z_3\}}\one{\{\z_3{\leftrightarrow}\z_4\}}$, the first factor is measurable with respect to $\F_{u_4-}$, while the conditional expectation of $\one{\{\z_3{\leftrightarrow}\z_4\}}$ is equal to $\phi(u_4 |y_3'|^d / f(Z_{\z_3}(u_4-)))$.
Using a first spatial integration with respect to $y_3'$, we get
$$\int dy_3' \,  \E[\one_{A_1} | \F_{u_4-}]= \one{\{\z_0{\leftrightarrow}\z_1{\leftrightarrow} \z_2{\leftrightarrow} \z_3\}} \, \frac {f(Z_{\z_3}(u_4-))} {u_4}. $$
The conditional expectation of the right hand side given $\F_{u_0}$ equals
\begin{align*}
\one{\{\z_0{\leftrightarrow}\z_1{\leftrightarrow} \z_2{\leftrightarrow} \z_3\}}  & \frac {(1+ Z_{\z_3}(u_0))} {u_4} \E\Big[\frac {f(Z_{\z_3}(u_4-))} {1+Z_{\z_3}(u_0)} \Big| \F_{u_0} \Big] \\
& \le c_1 \,  \one{\{\z_0{\leftrightarrow}\z_1{\leftrightarrow} \z_2{\leftrightarrow} \z_3\}} \frac {(1+ Z_{\z_3}(u_0))} {u_4}\Big (\frac {u_4} {u_0}\Big)^\gamma,
\end{align*}
where the  inequality follows from Lemma~\ref{boundonmoments} in the appendix. 
We now take conditional expectation given~$\F_{u_0-}$, note that $Z_{\z_3}(u_0)=Z_{\z_3}(u_0-)$,
and integrate in space with respect to $y_1'$, to obtain the bound
$$ c_1 \, \frac {f(Z_{\z_1}(u_0-))} {u_0} \one{\{\z_1{\leftrightarrow} \z_2{\leftrightarrow} \z_3\}} \frac {(1+ Z_{\z_3}(u_0-))} {u_4} \Big(\frac {u_4} {u_0}\Big)^\gamma.$$
The conditional expectation of this bound given $\F_{u_2}$ can be bounded by
$$ c_2 \, \frac {(1+Z_{\z_1}(u_2))} {u_0} \Big(\frac {u_0} {u_2}\Big)^\gamma \one{\{\z_1{\leftrightarrow} \z_2{\leftrightarrow} \z_3\}} \frac {(1+ Z_{\z_3}(u_2))} {u_4} \Big(\frac {u_4} {u_2}\Big)^\gamma,$$
by  Corollary~\ref{boundonmoments2} in the appendix.
Further, the conditional expectation of this expression given $\F_{u_2-}$, and integrated over $y_1$ and $y_3$, is exactly equal to
$$ c_2 \frac {(2+Z_{\z_1}(u_2-))} {u_0} \Big(\frac {u_0} {u_2}\Big)^\gamma \frac {f(Z_{\z_1}(u_2-))} {u_2} \frac {f(Z_{\z_3}(u_2-))} {u_2}
 \frac {(2+ Z_{\z_3}(u_2-))} {u_4} \Big(\frac {u_4} {u_2}\Big)^\gamma,$$
and bounded by
$$ c_3 \, \frac {(1+Z_{\z_1}(u_2-))} {u_0} \Big(\frac {u_0} {u_2}\Big)^\gamma \frac {(1+Z_{\z_1}(u_2-))} {u_2} \frac {(1+Z_{\z_3}(u_2-))} {u_2}
 \frac {(1+ Z_{\z_3}(u_2-))} {u_4} \Big(\frac {u_4} {u_2}\Big)^\gamma.$$
Using Corollary~\ref{boundonmoments2} again, we bound the conditional expectation given $\F_{u_1}$  
by
$$ c_4 \frac {1} {u_0} \Big(\frac {u_0} {u_1}\Big)^\gamma \frac {1} {u_2} \Big(\frac {u_2} {u_1}\Big)^\gamma \ \frac {(1+Z_{\z_3}(u_1))} {u_2} \Big(\frac {u_2} {u_1}\Big)^\gamma
 \frac {(1+ Z_{\z_3}(u_1))} {u_4} \Big(\frac {u_4} {u_1}\Big)^\gamma.$$
Finally, the expectation of that expression is bounded, using Lemma~\ref{boundonmoments} again,~by 
$$ c_5 \frac {1} {u_0} \Big(\frac {u_0} {u_1}\Big)^\gamma \frac {1} {u_2} \Big(\frac {u_2} {u_1}\Big)^\gamma \ \frac {1} {u_2} \Big(\frac {u_2} {u_3}\Big)^\gamma
 \frac {1} {u_4} \Big(\frac {u_4} {u_3}\Big)^\gamma.$$
Altogether,  we have proved that
$$\int_{(\R^d)^4} dz_1 dz_2 dz_3 dz_4 \,
\P (A_1)
\le c_5 \, \frac 1 {u_0^\gamma u_1^{1-\gamma}} \frac 1 {u_2^\gamma u_1^{1-\gamma}} \frac 1 {u_2^\gamma u_3^{1-\gamma}} \frac 1 {u_4^\gamma u_3^{1-\gamma}},$$
which gives the desired result if $C$ is chosen at least equal to $c_5^{1/4}$.
\medskip

In this calculation, space integration is used extensively to give a simple expression for the density of the probability, for a vertex $\z$ with indegree $k$ at time $t-$, to have a child somewhere with birth time in $dt$, 
$$\int_{\R^d} dz' \,  \P_{\z,(z',t)}\big(\z{\leftrightarrow}(z',t)\, \big|\, Z(\z,t-)=k\big) = f(k)/t.$$
It is important to perform the space integration at time $t$, before studying the indegree evolution 
process $Z(\z,s)$ for $s<t$. This method is independent of the choice of a profile function, showing 
that the argument does not involve space. But an alternative approach would be to use the profile function more explicitly. Because $\phi$ is regularly varying with index $\delta$, from the Potter bounds,  for any $\delta'<\delta$, there exists a finite constant $c$ such that
$\phi(x)\le c x^{-\delta'}$, for all $x>0$.  Fix a choice of such a $\delta'\in (1,\infty)$ and $c>0$.
Then, for given $\z=(z,u)$ older than $\z'=(z',u')$, we have
$$ \P_{\z,\z'}\big(\z{\leftrightarrow}\z'\ \big| \, Z(\z,u'-)=k \big)
\le c \,  f(k)^{\delta'} \left(u' \tmtextsf{d}(z,z')^{1/d}\right)^{-{\delta'}},$$
and then
$$ \P_{\z,\z'} \big(\z{\leftrightarrow}\z'\big)\le 1\wedge c\, \E[Z(\z,u'-)^{\delta'}] \left(u' \tmtextsf{d}(z,z')^{1/d}\right)^{-{\delta'}}.$$
With the same outline of proof, but using Lemma~\ref{boundonmoments} and Corollary~\ref{boundonmoments2} with 
$p={\delta'}$ or $p=2{\delta'}$, we could show that
\begin{equation} \label{quickpathbound_withoutspaceintegration}
\begin{aligned}
\P&_{(0,u_0),  (z_1,u_1),\ldots,(z_n,u_n)} \big( ((z_k,u_k))_{0\le k \le n} \text{ is a quick path in }G^\infty\big)\\
&\le c'^n \, \prod_{k=1}^n   1 \wedge \left((u_{k-1}\wedge u_k)^\gamma (u_{k-1}\vee u_k)^{1-\gamma} \vert z_k-z_{k-1}\vert\right)^{-{\delta'}},
\end{aligned}\end{equation}
for some constant~$c'$, and deduce Proposition~\ref{time_path_bound_prop} by integration over all the space variables. A similar bound will be used in the  next subsection, without further justification.

\subsection{Non--robustness for  $\delta>\frac1{1-\gamma}$  in dimension one}
We only need to consider $\gamma \in[1/2,1)$. In this phase, we always have $\delta>2$. 
We look for ways to improve the bound from the previous section. Any such argument
has to use the spatial structure of the network substantially, as the corresponding nonspatial 
networks are robust for $\gamma \geq 1/2$.
\medskip

Let us first sketch the idea informally. Suppose that $d=1$.
A vertex $\z=(z,u)$ has typically of order $u^{-\gamma}$ children, which may be a lot.
But most of these children are typically close to~$\z$, namely within distance $u^{-1}$, and hence their local neighbourhoods are 
strongly correlated.  No matter how many vertices within distance $u^{-1}$ of $\z$ belong to the component of $\z$, it will not help much to connect~$\z$ to vertices far away. Indeed, defining the \emph{region around $\z$} as
$$ C_\z:=\big\{(z',u'), u'\ge u, \vert z' -z \vert \le 2 u^{-1}-u'^{-1}\big\},$$
see Figure~4, we can show that the typical number of vertices outside $C_\z$ that are connected to $\z$, or {any other vertex} 
in $C_\z$, is only of order $\log(u^{-1})$. To estimate the probability of a path it therefore makes sense to consider 
only those edges of a quick path straddling the boundary of a region. This idea leads
us to the notion of a \emph{trace of a quick path} which we use to improve our bounds. Informally, forgetting about the time component and just thinking about $C_\z$ as a ball around $z\in \R$, it
is plausible that in dimension~$d=1$ few edges straddle the boundary of $C_\z$ because the size of the boundary of balls in $\R$ does not grow with the radius.   In dimension $d>1$ however, if we wanted to use a similar approach, we would have to consider the ball of radius $u^{-1/d}$ around a vertex. The area of its boundary is of order $u^{-(d-1)/d}$, and therefore the vertices within this ball would be connected to typically at least $O(u^{-(d-1)/d})$ vertices outside, which is already too much to carry out the proof.
\medskip

Suppose now that $(0,U)$ belongs to an infinite component of $\tensor[^p]{G}{^\infty}$.
Then it is the starting point of an infinite quick path $(\x_n)_{n\ge0}$ in $G^\infty$, as defined in
the previous section, in which every regular vertex is in $\tensor[^p]{\X}{}$.
We  define the subsequence $(\y_n)_{n\ge0}$ given by
\smash{$\y_n=\x_{\phi(n)}$} with $\phi(0)=0$, and 
$$\phi(n+1)= \min\{k>\phi(n),\x_k\notin C_{\x_{\phi(n)}}\}.$$
We call $(\y_n)_{0\leq n \leq m}$ a \emph{trace} of the quick path $(\x_n)_{n\ge0}$.
Observe that if $\y_n=\x_{\phi(n)}$ is a local maximum of the quick path, then $\y_{n+1}=\x_{\phi(n+1)}=\x_{\phi(n)+1}$ is regular and is in $\tensor[^p]{\X}{}$.
Therefore at least half of the vertices of the trace of a quick path are in $\tensor[^p]{\X}{}$.
Arguing in the same way as in last subsection, we have to prove that the expected number  
of traces of length~$n$ grows at most exponentially in~$n$.
This follows from the following two results, which are analogous to Proposition~\ref{time_path_bound_prop} and Lemma~\ref{time_integral_bound}.

\begin{proposition} \label{time_path_bound_prop2}
In dimension $d=1$, if $\frac12\leq\gamma<1$ and $\delta>\frac1{1-\gamma}$, then there exists a finite constant $C''$ such that, for every $n\ge0$ and $t_0,\ldots,t_n\in (0,1]$ pairwise distinct,  
\begin{equation} \label{time_path_bound2}
\begin{aligned}
\int dy_1 \ldots dy_n \,  &
\P_{\y_0, \ldots,\y_n} \big( (\y_k)_{0\le k \le n} \text{ is
a trace of a quick path in }G^\infty\big) \\
&\le C''^n \,  \prod_{k=1}^n \Big(\frac {\one\{t_k<t_{k-1}\}} {t_{k-1}^{1-\gamma} t_k^{\gamma}} +
 \frac {\one\{t_k>t_{k-1}\}} {t_{k}}\Big),
\end{aligned}
\end{equation}
where we have written $\y_0=(0,t_0)$ and $\y_k=(y_k,t_k)$, and where the domain of integration is
$\{(y_1,\ldots,y_n)\in \R^n \colon$  $\y_k \notin C_{\y_{k-1}}$ for all $k \}$.
\end{proposition}

The improvement of this bound, compared to~\eqref{time_path_bound}, is that if 
$t_{k-1}<t_k$, the term \smash{$1/t_{k-1}^{\gamma} t_k^{1-\gamma}$} has been replaced by $1/t_k$.
Note also that the proposition is valid only in dimension one and for parameters satisfying $\gamma \in[1/2,1)$ and $\delta>1/(1-\gamma)$.

\begin{lemma} \label{time_integral_bound2}
There exists a finite constant $C'''$ such that, for any $n>0$, we have
$$ \int_0^1 dt_0 \cdots \int_0^1 dt_n \,  \prod_{k=1}^n
\Big(\frac {\one\{t_k<t_{k-1}\}} {t_{k-1}^{1-\gamma} t_k^{\gamma}} +
 \frac {\one\{t_k>t_{k-1}\}} {t_{k}}\Big) \le C'''^n.$$
\end{lemma}
\medskip

We skip the simple proof of Lemma~\ref{time_integral_bound2}. To prove Proposition~\ref{time_path_bound_prop2}, we  do not  bound the probability of a sequence being the the trace of a quick path $(\x_k)_{k\ge 0}$ directly, but instead construct a third sequence~\smash{$(\z_k)_{k\ge 0}$} called an \emph{almost quick path}.   To this end we first
define  the \emph{enlarged region} \smash{$C'_{\y_{k-1}}$} around vertex \smash{$\y_{k-1}=(y_{k-1},t_{k-1})$} by
\begin{align*}
C'_{\y_{k-1}}:=\big\{ & (y,t)\in \R\times(0,1] \colon t\ge t_{k-1}, \tmtextsf{d}(y,y_{k-1})\le 2 t_{k-1}^{-1}+t^{-1}\big\}  \\
&\cup \, \big\{(y,t)\in \R\times(0,1] \colon  t< t_{k-1}, \tmtextsf{d}(y,y_{k-1})\le 2 t_{k-1}^{-1}+t_{k-1}^{\gamma-1} t^{-\gamma}\big\}.\\[-10mm]
\end{align*}
\begin{figure}[ht]
\begin{center}
\includegraphics[scale=0.6]{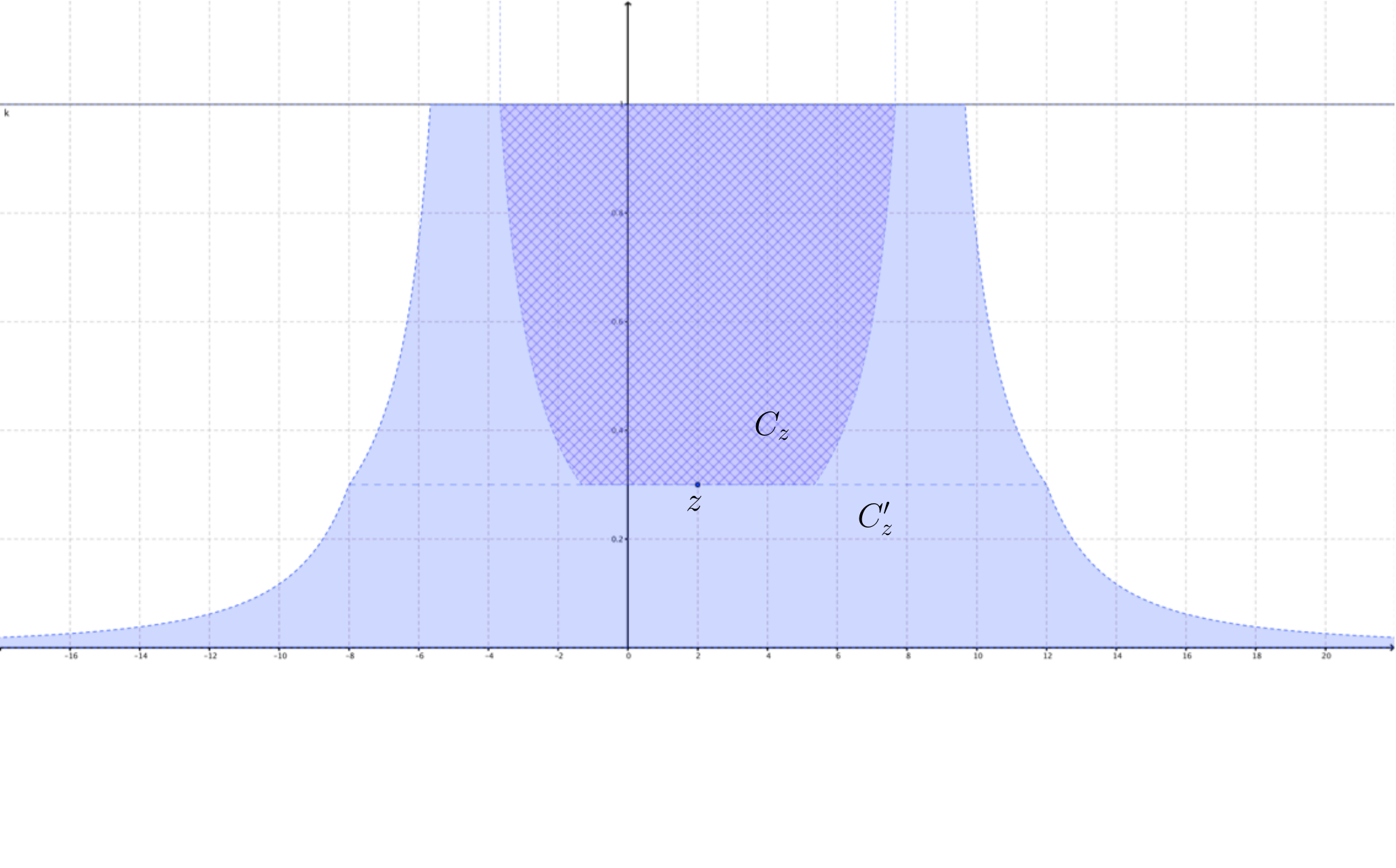}
\end{center}
\ \vspace{-12mm}\\ 
\caption{Graphical representation of the region $C_\z$ (dark shade) and enlarged region $C'_\z$ (light shade) around the  vertex $\z=(2,0.3)\in \R\times(0,1]$, for~$\gamma=1/2$.\vspace{-3mm}}
\end{figure}

Define $(\z_k)_{k\ge 0}$ by inserting in the infinite trace $(\y_k)_{k\ge 0}$ a vertex \smash{$\y'_k:=\x_{\phi(k)-1}$} between the vertices \smash{$\y_{k-1}=\x_{\phi(k-1)}$} and $\y_{k}=\x_{\phi(k)}$, but only if 
\begin{itemize}
\item $\phi(k)-1>\phi(k-1)$ and 
\item $\y_{k} \notin C'_{\y_{k-1}}$. 
\end{itemize}
In other words, if  $\y_k$ is even outside the enlarged region~\smash{$C'_{\y_{k-1}}$} and it is not already represented, we insert the vertex in $C_{\y_{k-1}}$ connecting to $\y_{k}$.
The infinite sequence $(\z_k)$ we obtain is again a subsequence of the quick path \smash{$(\x_k)$}.
Again a vertex $\z_k$ is called regular for this sequence if it is older than $\z_{k+1}$ or $\z_{k-1}$, and otherwise it is called a local maximum.
Observe that local maxima of the sequence $(\z_k)$ are not necessarily local maxima of the sequence $(\x_k)$, but regular vertices of $(\z_k)$ are always regular vertices of $(\x_k)$.
\medskip

\pagebreak[3]

It is not hard to show that the sequence $(\z_k)$ satisfies Properties (ii) and (iii) of the definition of quick paths.
Actually, it can fail to be a quick path itself, only because it may not be a path, as some of the pairs $(\z_k,\z_{k+1})$ are not requested to be edges of the graph. Note though that from the sequence $(\z_k)$, one can identify the vertices $\y_k$, the inserted vertices $\y'_k$, and which pairs $(\z_k,\z_{k+1})$ are required to be edges of the graph, and which are not. A self-avoiding sequence satisfying Properties (ii), and (iii), such that 
all requested edges are present is called an \emph{almost quick path}.

\subsubsection*{Proof of Proposition~\ref{time_path_bound_prop2}}
Fix $n\ge0$ and $t_0,\ldots,t_n$ distinct times in $(0,1]$.
Let $y_1,\ldots,y_n$ be real numbers such that, defining $\y_0=(0,t_0)$ and $\y_k=(y_k,t_k)$ for $1\le k \le n$, the sequence $(\y_k)_{0\le k\le n}$ satisfies $\y_k \notin C_{y_{k-1}}$ for all $1\le k \le n$.
Let 
$$A=\big\{ k\in\{1,\ldots,n\} \colon \y_k\notin C'_{y_{k-1}}\big\}.$$
If the sequence $(\y_k)_{k\le n}$ is the trace of a quick path, there must be $B\subset A$, say of cardinality~$m$,  
and, for every $k\in B$, a vertex 
$\y'_k \in C_{y_{k-1}}$, such that the sequence $(\z_k)_{0\le k \le n+m}$ obtained by inserting the vertices $\y'_k$, 
is an almost quick path.
Consequently, we have
\begin{align*}
\P_{\y_0, \ldots,\y_n}  & \big( (\y_k)_{0\le k \le n} \text{ is
the trace of a quick path in }G^\infty\big) \\
&\le \sum_{B\subset A} \int d\y'_{k_1} \ldots d\y'_{k_m}
\P_{\z_0,\ldots,\z_{n+m}}\big((\z_k)_{0\le k \le n+m} \text{ is an almost quick path}\big),
\end{align*}
where we have written $k_1,\ldots, k_m$ for the ordered elements of $B$.
The number of pairs $(A,B)$ with $B\subset A \subset \{1,\ldots,n\}$ is equal to $3^{n}$, thus in order 
to prove~\eqref{time_path_bound2} it suffices to show that, for any possible choice of $A$ and $B=\{k_1,\ldots, k_m\}\subset A$, 
we have
\begin{align*}
\int dy_1 \ldots dy_n & \, \int d\y'_{k_1} \ldots d\y'_{k_m}
\P_{\z_0,\ldots,\z_{n+m}}\big((\z_k)_{0\le k \le n+m} \text{ is an almost quick path} \big)\\
&\le \Big( \frac {C''} 3\Big)^n \prod_{k=1}^n \Big(\frac {\one\{t_k<t_{k-1}\}} {t_{k-1}^{1-\gamma} t_k^{\gamma}} +
 \frac {\one\{t_k>t_{k-1}\}} {t_{k}}\Big),
\end{align*}
where the domain of integration, depending on $A$ and $B$, is defined by the constraints
\smash{$\y_k\in C'_{\y_{k-1}}\backslash C_{\y_{k-1}}$ for $k\notin A$, $\y_k \notin C'_{\y_{k-1}}$
for $k\in A$, and $\y'_k \in C_{\y_{k-1}}$ for $k\in B$.}
\medskip

\pagebreak[3]

We first give a bound on the probability that $(\z_k)_{0\le k \le n+m}$ is an almost quick path in the same way as for a quick path.
We keep now the notation $(z_k,u_k)$, resp. $(y_k,t_k)$ and $(y'_k,t'_k)$, for $\z_k$, resp. $\y_k$ and $\y'_k$, for appropriate indices $k$.
For $\frac1{1-\gamma}<\delta'<\delta$ we replace \eqref{quickpathbound_withoutspaceintegration} by
\begin{equation} \label{almostquickpathbound}
\begin{aligned}
\P_{\z_0,\ldots,\z_{n+m}} & \big((\z_k)_{0\le k \le n+m} \text{ is an almost quick path} \big)\\
&\le 
c^n \prod_{k\in B} \left((t'_k\wedge t_k)^\gamma (t'_k\vee t_k)^{1-\gamma} \vert y_k-y'_k\vert\right)^{-{\delta'}}\\
& \phantom{ganzlangewarten}\times \prod_{k\in A\backslash B} \left((t_{k-1}\wedge t_k)^\gamma (t_{k-1}\vee t_k)^{1-\gamma} \vert y_k-y_{k-1}\vert\right)^{-{\delta'}},
\end{aligned}\end{equation}
for $c=c(\delta')$ some finite constant.  Observe that each factor corresponds to a requested edge. 
The proof of~\eqref{almostquickpathbound} requires to split the paths into small parts and then give a bound 
for the individual probability of each part.  We do not provide the detail of this proof, as this is very similar to the previous section.
Instead, we now show how to perform the integration over the variables $y_k$ and $\y'_k$ to get an improved bound.
\medskip

\pagebreak[3]

We first introduce the change of variables \smash{$\tilde y_k=y_k-y_{k-1}$} for $1\le k\le n$ and \smash{$\tilde y'_k=y'_k-y_{k-1}$} for $k\in B$, and write \smash{$\tilde \y_k=(\tilde y_k,t_k)$ and $\tilde \y'_k=(\tilde y'_k,t'_k)$}. The domain of integration is now defined by the constraints
\smash{$\tilde \y_k\in C'_{(0,t_{k-1})}\backslash C_{(0,t_{k-1})}$} for $k\notin A$, \smash{$\tilde \y_k \notin C'_{(0,t_{k-1})}$} for $k\in A$, and
\smash{$\tilde \y'_k \in C_{(0,t_{k-1})}$} for $k\in B$,
which is a product domain with respect to the variables $\tilde y_k$ and $\tilde \y'_k$.
Proposition~\ref{time_path_bound_prop2} will follow if, for each $k$, we can integrate over $\tilde y_k$ (resp.\ over $(\tilde \y'_k, \tilde y_k)$ if $k\in B$), the single term in the product on the right hand side of~\eqref{almostquickpathbound} involving this variable, and ensure the result is bounded by a constant multiple~of
$$ \frac {\one\{t_k<t_{k-1}\}} {t_{k-1}^{1-\gamma} t_k^{\gamma}} + \frac {\one\{t_k>t_{k-1}\}} {t_{k}}.$$
This is what we now do, considering separately the different cases.
\medskip

{\bf (A) The case~$k\notin A$.} 
Integrating a constant over  the domain \smash{$\{\tilde y_k \colon \tilde \y_k\in C'_{(0,t_{k-1})} \backslash C_{(0,t_{k-1})}\}$}
we obtain a term of order \smash{$1/t_{k}$}, if $t_{k-1}<t_k$, and of order \smash{$1/t_{k-1}^{1-\gamma} t_k^{\gamma}$}, if $t_{k-1}>t_k$.
\smallskip

{\bf (B) The case $k\in A\backslash B$.} Then we have to integrate
\smash{$((t_{k-1}\wedge t_k)^\gamma (t_{k-1}\vee t_k)^{1-\gamma} \vert \tilde y_k\vert)^{-{\delta'}}$}
over the domain \smash{$\{\tilde y_k \colon \tilde \y_k\notin C'_{(0,t_{k-1})}\}$}. The reader can easily check
the bound in this case.
\smallskip


{\bf (C1) The case $k\in B$ and $t_{k-1}>t_k$.} We have to integrate 
\smash{$t_k^{-\gamma\delta'} (t'_k)^{-(1-\gamma) \delta'} \vert \tilde y_k - \tilde y'_k\vert^{-\delta'}$}
over \smash{$t'_k>t_{k-1}$}, \smash{$\vert \tilde y'_k\vert \le 2 t_{k-1}^{-1}-(t'_k)^{-1}$}, and 
\smash{$\vert \tilde y_k\vert \ge 2 t_{k-1}^{-1}+ t_{k-1}^{\gamma-1}t_k^{-\gamma}$}.
Write $u:=\vert\tilde y_k\vert -2 t_{k-1}^{-1}$ and bound  \smash{$\vert \tilde y_k - \tilde y'_k\vert^{-\delta'}$ by $u^{-\delta'}$}, so that the integral over $\tilde y'_k$ gives at most a factor $4t_{k-1}^{-1}$, and the integral is bounded by
$$8 t_k^{-\gamma\delta'} (t_{k-1})^{-1} \int_{t_{k-1}}^{1} (t'_k)^{-(1-\gamma) \delta'} dt'_k \int_{t_{k-1}^{\gamma-1}t_k^{-\gamma}}^{\infty} u^{-\delta'} du\le c t_{k-1}^{\gamma-1} t_k^{-\gamma},$$
for some finite constant $c$. We have used that $\delta'>1$ and $(1-\gamma)\delta'>1$ to obtain the right order for the integrals in~\smash{$t'_k$} and in~$u$.
\smallskip

{\bf (C2) The case $k\in B$ and $t_{k-1}<t_k$.} 

First, we bound the integral over $\tilde \y'_k$ younger than $\tilde\y_k$.
We have to integrate again the quantity \smash{$t_k^{-\gamma\delta'} (t'_k)^{-(1-\gamma) \delta'} \vert \tilde y_k - \tilde y'_k\vert^{-\delta'},$}
but the integration is now over \smash{$t'_k>t_{k}$, $\vert \tilde y'_k\vert \le 2 t_{k-1}^{-1}-(t'_k)^{-1}$}, and
\smash{$\vert \tilde y_k\vert \ge 2 t_{k-1}^{-1}+ t_k^{-1}$.}
Writing \smash{$u:=\vert\tilde y_k\vert -2 t_{k-1}^{-1}$} and 
$v:=\vert \tilde y_k-\tilde y'_k\vert$, we can similarly bound the integral by
$$2 t_k^{-\gamma\delta'} \int_{t_{k}}^1 (t'_k)^{-(1-\gamma) \delta'} dt'_k
\int_{t_k^{-1}}^{\infty} \left(\int_{u}^{\infty} v^{-\delta'} dv\right) du
\le c t_k^{1- \delta'} \int_{t_k^{-1}}^{\infty} u^{1-\delta'} du  \le c' t_k^{-1},$$
for $c,c'$ some finite constants. We have used the fact $(1-\gamma)\delta'>1$ to obtain the order of the integral in $t'_k$, 
and we have used $\delta'>2$ to bound the integral in $u$.
\smallskip

Second, we bound the integral over $\tilde \y'_k$ older than $\tilde \y_k$.
The quantity we have to integrate is now
\smash{$(t'_k)^{-\gamma\delta'} t_k^{-(1-\gamma) \delta'} \vert \tilde y_k - \tilde y'_k\vert^{-\delta'},$}
and the integration is over \smash{$t_{k-1}<t'_k<t_{k}$, $\vert \tilde y'_k\vert \le 2 t_{k-1}^{-1}-(t'_k)^{-1}$}, and
\smash{$\vert \tilde y_k\vert \ge 2 t_{k-1}^{-1}+ t_k^{-1}.$}
With the same notation as before we have \smash{$v>(t'_k)^{-1}$ and $u<v$}, and the integral in $u$, for $v$ fixed, 
gives at most a factor $v$, so that the integral is bounded by
\begin{align*}
2 t_k^{-(1-\gamma)\delta'} \! \int_{t_{k-1}}^{t_k}  (t'_k)^{-\gamma \delta'}\Big(
\int_{(t'_{k})^{-1}}^{\infty} v^{1-\delta'} dv\Big) dt'_k
\le  c t_k^{-(1-\gamma)\delta'}
\int_{t_{k-1}}^{t_k} (t'_k)^{(1-\gamma) \delta'-2} dt'_k
\le c' t_k^{-1},
\end{align*}
for some constants $c, c'$.
We have used that $1-\delta'<-1$ to obtain the right order for the integral in $v$, and that $(1-\gamma)\delta'-2>-1$ to obtain the right order for the integral in $t'_k$.

\section{Appendix: Auxiliary lemmas} \label{auxiliaries}

For each $t\in [1,\infty]$ fixed, the graph $G^t$ is constructed as a growing graph with vertices placed in $\T_t$ and with birth times in $(0,1]$. The indegree of a vertex $\x=(x,r)$ at time $s\ge r$ is denoted by $Z^t(\x, s)$. The process
$(Z^t(\x, s))_{s\ge r}$ is  a time-inhomogeneous pure birth process started in zero at time $s=r$. By translation in 
$\T_t$, the law of this process does not depend on $x$, and we write $Z^t(r, s)$ for $Z^t((0,r), s)$ under 
the measure~$\P_{(0,r)}$, i.e., conditionally on the vertex we consider to be in the Poisson point process.
\smallskip

This appendix provides different estimates and bounds for this process. We treat simultaneously the cases $t=1$ and $t>1$, including $t=\infty$, and do not stop the process at time $s=1$. 
The process $(Z^t(r, s))_{s\ge r}$ was already studied in~\cite{JM14}, see in particular Lemma~8, where we proved 
that $\log Z^t(r,s)\sim \gamma \log s$ almost surely as $s\uparrow \infty$. Moreover, it 
was shown\footnote{This statement is actually proved under a slightly stronger assumption on that attachment rule~$f$, namely that
$f(k)=\gamma k + O(1)$. As explained in \cite[Remark 6]{JM14} the results carry over to our framework at the price of an arbitrarily small 
increase of $\gamma$. This is still sufficient for our applications.}   in Lemma~9 that the probability of having a larger indegree than $(s/r)^\gamma$ decays exponentially, namely
\begin{equation}\label{upperboundindegree}
\P(Z^t(r,s)\ge \lambda (s/r)^\gamma) \le c \exp(-\lambda/8),
\end{equation}
for some (explicit) constant $c$ only depending on the attachment rule $f$. An easy modification of the argument gives a similar result for the increase of the process on the interval~$(s,s']$,~i.e.\
\begin{equation}\label{upperboundindegreeincrease}
\P\left(\frac {1+Z^t(r,s')}{1+Z^t(r,s)} \ge \lambda (s'/s)^\gamma \Big| \, Z^t(r,s)\right) \le c \exp(-\lambda/8).
\end{equation}
Here it is important that the bound does not depend on the value taken by $Z^t(r,s)$. A consequence of this exponentially decaying tail is that the moments are well-controlled, see the next lemma and its corollary.
\begin{lemma}\label{boundonmoments}
For each $p\in [1,\infty)$, 
there exists a constant $c_p$ depending only on $p$ and on the attachment rule, so that for every $r\le s<s'$, we have
$$\E\Big[\Big( \frac {1+Z^t(r,s')}{1+Z^t(r,s)}\Big)^p \Big|Z^t(r,s)\Big] \le c_p \Big(\frac {s'} s \Big)^{p\gamma}.$$
\end{lemma}
\begin{proof}
For a positive random variable $X$, we have
$$ \E[X^p]\le \sum_{k=0}^\infty \P(k\le X<k+1) (k+1)^p \le\sum_{k=0}^\infty \P(X\ge k) (k+1)^p,$$
which is bounded by an explicit finite constant if $X$ has an explicit exponentially decaying bound. Apply this to the random variable
$$ \frac {1+Z^t(r,s')}{1+Z^t(r,s)} \left(\frac s {s'}\right)^\gamma.$$
\end{proof}
\begin{corollary}\label{boundonmoments2}
Let $\x$ and $\y$ be vertices born before time $s$, and $s'>s$. For each~$p\in [1,\infty)$, 
$$\E\left[\left( \frac {1+Z^t(\x,s')}{1+Z^t(\x,s)}\right)^p \left( \frac {1+Z^t(\y,s')}{1+Z^t(\y,s)}\right)^p 
\Big| \, Z^t(\x,s), Z^t(\y,s)\right] \le c_{2p} \left(\frac {s'} s \right)^{2p\gamma},$$
where $c_{2p}$ is as introduced in Lemma~\ref{boundonmoments}.
\end{corollary}
\begin{proof}
Use the Cauchy-Schwarz inequality, followed by Lemma~\ref{boundonmoments} .
\end{proof}

The next lemma gives a bound on the probability of observing a \emph{small} degree.

\begin{lemma} \label{lower_bound} There exists a function $\tilde g\colon(0,\infty)\to(0,\infty)$ growing at infinity slower than polynomially, 
such that
\begin{equation} \label{lower_bound_indegree}
 \P\big(Z^1(1,r)\le r^\gamma/ \tilde g(r)\big) \underset {r\to \infty}{\longrightarrow} 0.
  \end{equation}
\end{lemma}

\begin{proof}
It was proved in~\cite{JM14} that $\log(Z^{1}(1,r))/(\gamma \log r)$ converges almost surely and in probability to one. 
In particular, there exists a function $\psi$ such that for any $\eta>0$ and any $r\ge \psi(\eta)$, we have
$\P(Z^{1}(1,r)\le r^{\gamma (1-\eta)}) \le \eta$.
The function $\psi$ can be chosen decreasing with infinite limit at zero, so that its inverse~$\psi^{-1}$ is decreasing 
and converging to zero. For any $r>1$, we thus have
$\P(Z^{1}(1,r)\le r^{\gamma} r^{-\gamma \psi^{-1} (r)}) \le \psi^{-1}(r).$
Hence we can choose $\tilde g(r)= r^{\gamma \psi^{-1} (r)}$, which is 
$o(r^l)$ for any~$l>0$.
\end{proof}

\begin{lemma} \label{goodness_prob_bound} 
There exists a function $g\colon(0,\infty)\to(0,\infty)$ growing at infinity slower than polynomially, such that
\begin{equation} \label{good_vertex_inequality}
 \sup_{\heap{t \in (1,+\infty]}{\frac 1 {t \log t}\le s}}
\P\big(Z^t(s,\sfrac12)\le s^{-\gamma}/g(s^{-1})\big) \underset {s\to 0}{\longrightarrow} 0.
  \end{equation}
\end{lemma}

\begin{proof}
The supremum over $t$ is attained when $t$ is smallest possible, i.e.\ $t\log t=\frac1s$. Using that $Z^t(r',\frac12)$ is stochastically dominated by $Z^t(r,\frac12)$ if $r<r'$, we have, for $t>e$,
\begin{eqnarray*}
\P\big(Z^t(\sfrac1{t\log t}, \sfrac12)\le (t\log t)^{\gamma}/g(t\log t)\big)
&\le&\P\big(Z^t(\sfrac1t, \sfrac12)\le (t\log t)^{\gamma}/g(t\log t)\big)\\
&\le& \P\big(Z^1(1,\sfrac{t}2)\le (t\log t)^{\gamma}/g(t\log t)\big),
\end{eqnarray*}
where the second line follows from the scaling property. In order to prove the result, using~\eqref{lower_bound_indegree}, it is enough to ensure that we can choose $g$ growing slower than polynomially such that
$(t \log t)^\gamma / g(t \log t) \le (t/2)^\gamma/\tilde g(t/2)$, e.g.\ by letting
$g(u)=(2 \log u)^\gamma \sup_{t\le u/2} \tilde g(u).$
\end{proof}

We stress that the probability of having a smaller degree than expected does not decay exponentially. 
Indeed, the probability that $Z^\infty(s,1)=0$, i.e.\ the indegree of the vertex born at time~$s$ is still null at time $1$, 
decays only polynomially in $s^{-1}$.
Hence, despite the fact that a vertex born at time~$s$ typically has total indegree $s^{- \gamma+o(1)}$, 
there may well be some untypical vertices with much fewer inedges.
\bigskip

{\bf Acknowledgements:} This work was initiated when EJ visited Bath funded by a grant from 
the \emph{ESF} programme on `Random Geometry of Large Interacting Systems and Statistical Physics' (RGLIS) in November~2012.
PM was supported by \emph{EPSRC} grant~EP/K016075/1,  and by an invitation to  ENS Lyon in September~2013. 
The work  was completed when EJ spent the last four months of 2014 at the University  of Bath, funded by~\emph{CNRS}. We would like to 
thank Rob van den Berg and G\"unter Last for useful discussions
\medskip

\bibliographystyle{abbrv}
\bibliography{spatial_rob_biblio}

\end{document}